\documentclass[11pt]{article}

\usepackage{amsmath, amsfonts, amssymb,latexsym}

\usepackage{amsthm}
\usepackage{bbm}

\input epsf.sty

\DeclareFontFamily{U}{rcjhbltx}{}
\DeclareFontShape{U}{rcjhbltx}{m}{n}{<->rcjhbltx}{}
\DeclareSymbolFont{hebrewletters}{U}{rcjhbltx}{m}{n}

\let\aleph\relax\let\beth\relax
\let\gimel\relax\let\daleth\relax

\DeclareMathSymbol{\aleph}{\mathord}{hebrewletters}{39}
\DeclareMathSymbol{\beth}{\mathord}{hebrewletters}{98}
\DeclareMathSymbol{\gimel}{\mathord}{hebrewletters}{103}
\DeclareMathSymbol{\daleth}{\mathord}{hebrewletters}{100}
\DeclareMathSymbol{\he}{\mathord}{hebrewletters}{104}
\DeclareMathSymbol{\zayin}{\mathord}{hebrewletters}{122}
\DeclareMathSymbol{\kaf}{\mathord}{hebrewletters}{107}
\DeclareMathSymbol{\lamed}{\mathord}{hebrewletters}{108}
\DeclareMathSymbol{\mem}{\mathord}{hebrewletters}{109}
\DeclareMathSymbol{\noun}{\mathord}{hebrewletters}{110}
\DeclareMathSymbol{\samekh}{\mathord}{hebrewletters}{115}
\DeclareMathSymbol{\ayin}{\mathord}{hebrewletters}{96}
\DeclareMathSymbol{\pe}{\mathord}{hebrewletters}{112}
\DeclareMathSymbol{\tsade}{\mathord}{hebrewletters}{118}
\DeclareMathSymbol{\qof}{\mathord}{hebrewletters}{113}
\DeclareMathSymbol{\resh}{\mathord}{hebrewletters}{114}
\DeclareMathSymbol{\shin}{\mathord}{hebrewletters}{152}
\DeclareMathSymbol{\tav}{\mathord}{hebrewletters}{116}

\newtheorem{Theorem}{Theorem}[section]
\newtheorem{Lemma}{Lemma}[section]
\newtheorem{Proposition}[Lemma]{Proposition}
\newtheorem{Corollary}[Lemma]{Corollary}
\newtheorem{Definition}[Lemma]{Definition}

\newcommand{\BEQ}{\begin{equation}}     
\newcommand{\BEA}{\begin{eqnarray}}
\newcommand{\BD}{\begin{displaymath}}
\newcommand{\EEQ}{\end{equation}}       
\newcommand{\EEA}{\end{eqnarray}}
\newcommand{\ED}{\end{displaymath}}

\newcommand{\del}{\delta}
\newcommand{\Del}{\Delta}
\newcommand{\eps}{\varepsilon}          

\newcommand{\supp}{{\mathrm{supp}}}

\newcommand{\Tr}{{\mathrm{Tr}}}

\newcommand{\sh}{{\mathrm{sh}}}
\newcommand{\ch}{{\mathrm{ch}}}
\newcommand{\R}{\mathbb{R}}
\newcommand{\C}{\mathbb{C}}
\newcommand{\Z}{\mathbb{Z}}

\def\proba{{\mathbb{P}}}

\def\esper{{\mathbb{E}}}

\def\Cov{{\mathrm{Cov}}}

\def\sgn{{\mathrm{sgn}}}

\def\cotan{{\mathrm{cotan}}}

\newcommand{\ccotan}{{\widetilde{\cot}}}

\newcommand{\eop}{\hfill $\Box$}        

\newcommand{\II}{{\rm i}}               
\renewcommand{\Re}{{\rm Re\ }}          
\renewcommand{\Im}{{\rm Im\ }}          
\newcommand{\half}{{1\over 2}}          

\catcode`\@=11
\def\numberbysection{\@addtoreset{equation}{section}
        \def\theequation{\thesection.\arabic{equation}}}
\numberbysection

\begin{document}

\vspace*{1.5cm}
\begin{center}
{\Large \bf Global fluctuations  for 1D log-gas dynamics. (2) Covariance kernel and
support}

\end{center}


\vspace{2mm}
\begin{center}
{\bf  J\'er\'emie Unterberger$^a$}
\end{center}

\vskip 0.5 cm
\centerline {$^a$Institut Elie Cartan,\footnote{Laboratoire 
associ\'e au CNRS UMR 7502} Universit\'e de Lorraine,} 
\centerline{ B.P. 239, 
F -- 54506 Vand{\oe}uvre-l\`es-Nancy Cedex, France}
\centerline{jeremie.unterberger@univ-lorraine.fr}

\vspace{2mm}
\begin{quote}

\renewcommand{\baselinestretch}{1.0}
\footnotesize
{ 

We consider the hydrodynamic limit in the macroscopic regime of the coupled system of stochastic
differential equations,
\BEQ d\lambda_t^i=\frac{1}{\sqrt{N}} dW_t^i - V'(\lambda_t^i) dt+ \frac{\beta}{2N}
\sum_{j\not=i} \frac{dt}{\lambda^i_t-\lambda^j_t}, \qquad i=1,\ldots,N, \EEQ
with $\beta>1$, sometimes called {\em generalized Dyson's Brownian motion}, describing the dissipative
dynamics of a log-gas of $N$ equal charges with equilibrium measure corresponding to a $\beta$-ensemble, with
sufficiently regular convex potential $V$. The limit $N\to\infty$ is known to satisfy a mean-field
Mc Kean-Vlasov equation. Fluctuations around this limit have been shown \cite{Unt1} to
define a Gaussian process solving some explicit martingale problem written in terms of a generalized
transport equation.

We prove a series of results concerning either the Mc Kean-Vlasov equation for the
density $\rho_t$, notably regularity results and time-evolution of the support, or the associated hydrodynamic fluctuation process, whose 
space-time covariance kernel we compute explicitly.
 }

\end{quote}

\vspace{4mm}
\noindent

 \medskip
 \noindent {\bf Keywords:} random matrices, Dyson's Brownian motion, log-gas, beta-ensembles, hydrodynamic limit, Stieltjes transform,
 fluctuations, support

\smallskip

\noindent
{\bf Mathematics Subject Classification (2010):}  60B20; 60F05; 60G20; 60J60; 60J75; 60K35

\newpage

\tableofcontents



\section{Introduction and statement of main results}


\subsection{Introduction}


Let $\beta\ge 1 $ be a fixed parameter, and  $N\ge 1$ an integer. We consider the following system of coupled stochastic differential equations
driven by $N$ independent standard Brownian motions $(W_t^1,\ldots,W_t^N)_{t\ge 0}$,

\BEQ d\lambda_t^i=\frac{1}{\sqrt{N}} dW_t^i - V'(\lambda_t^i) dt+ \frac{\beta}{2N}
\sum_{j\not=i} \frac{dt}{\lambda^i_t-\lambda^j_t}, \qquad i=1,\ldots,N  \label{eq:SDE} \EEQ

Letting
\BEQ {\cal W}(\{\lambda^i\}_i):=\sum_{i=1}^N V(\lambda^i)-\frac{\beta}{4N} \sum_{i\not=j}
\log(\lambda^i-\lambda^j),\EEQ
we can rewrite (\ref{eq:SDE}) as $d\lambda_t^i=\frac{1}{\sqrt{N}} dW_t^i-\nabla_i {\cal W}
(\lambda_t^1,\ldots,\lambda_t^N) dt$.
Thus the corresponding equilibrium measure,
\BEQ d\mu^N_{eq}(\{\lambda^i\}_i) = \frac{1}{Z^N_V} e^{-2N{\cal W}(\{\lambda^i\}_i)}=\frac{1}{Z^N_V} \left(\prod_{j\not=i}  |\lambda^j-\lambda^i|\right)^{\beta/2} \exp\left(-2N\sum_{i=1}^N V(\lambda^i)\right)
\ d\lambda^1\cdots d\lambda^N \label{eq:eq-measure} \EEQ
is that of a $\beta$-log gas with confining potential $V$. 

\medskip

\noindent Let us start with a historical overview of the subject as a motivation for our study. This system of equations was originally considered in a particular case by Dyson \cite{Dys} who wanted to
describe the Markov evolution of a Hermitian matrix $M_t$ with i.i.d. increments $dG_t$ taken from the 
Gaussian unitary ensemble (GUE). In Dyson's idea, this matrix-valued process was to be
a matrix analogue of Brownian motion. The latter time-evolution being invariant through conjugation
by unitary matrices, we may project it onto a time-evolution of the set of eigenvalues 
$\{\lambda_t^1,\ldots,\lambda_t^N\}$ of the matrix, and obtain (\ref{eq:SDE}) with
$\beta=2$ and $V\equiv 0$. Keeping $\beta=2$, it is easy to prove that (\ref{eq:SDE}) 
is  equivalent to a generalized matrix Markov evolution, $dM_t=dG_t-V'(M_t) dt$.  The Gibbs measure
$${\cal P}^N_V(M)=\frac{1}{Z_N} e^{-N\Tr V(M)} dM, \qquad dM=\prod_{i=1}^N dM_{ii}
\prod_{1\le i<j\le n} d\Re M_{ij} \ d\Im M_{ij} $$
 can then be proved to be an  equilibrium measure. Such measures, together with their
projection onto the eigenvalue set, $\mu^N_{eq}(\{\lambda^1,\ldots,\lambda^N\})$, are the main
object of random matrix theory, see e.g. \cite{Meh},\cite{And}, \cite{Pas}. The {\em equilibrium
eigenvalue distribution} can be studied by various means, in particular using orthogonal
polynomials with respect to the weight $e^{-NV(\lambda)}$. The scaling in $N$ (called
{\em macroscopic scaling} in random matrix theory) ensures
the convergence of the random point measure $X^N:=\frac{1}{N} \sum_{i=1}^N \del_{\lambda^i}$ to  a deterministic  measure $\mu_V$ with {\em compact} support and density $\rho$ when $N\to\infty$ (see
e.g. \cite{Joh}, Theorem 2.1). One finds e.g. the well-known semi-circle law, $\rho(x)=\frac{1}{\pi}\sqrt{2-x^2}$, when $V(x)=x^2/2$. Looking more closely at the limit of the point
measure, one finds for arbitrary {\em polynomial} $V$  (Johansson \cite{Joh}) Gaussian fluctuations of order $O(1/N)$, contrasting with the
$O(1/\sqrt{N})$ scaling of fluctuations for the means of $N$ independent random variables,
typical of the central limit theorem. Assuming that the support of the measure is
connected (this essential "one-cut" condition holding in particular for $V$ {\em convex}), Johansson proves that the {\em covariance} of the limiting law depends on $V$ only through the support of the measure -- it is thus
{\em universal} up to a scaling coefficient --, while the means is
equal to $\rho$, plus an apparently non-universal correction in $O(1/N)$. 

\medskip

Following Rogers and Shi \cite{RogShi}, Li, Li and Xie \cite{LiLiXie} proved
the following two facts:
\begin{itemize}
\item[(i)]  two arbitrary eigenvalues never collide, which implies the non-explosion of (\ref{eq:SDE});

\item[(ii)] the random point process $X^N_t:=\frac{1}{N} \sum_{i=1}^N \del_{\lambda^i_t}$
satisfies in the limit $N\to\infty$ a deterministic hydrodynamic equation of Mc Kean Vlasov type, namely, the asymptotic density 
\BEQ \rho_t \equiv X_t:={\mathrm{w\!-\!\!\lim}}_{N\to\infty} X_t^N \EEQ
 satisfies the PDE
\BEQ \frac{\partial \rho_t(x)}{\partial t}=\frac{\partial}{\partial x}\left( \left( V'(x)-
\frac{\beta}{2} p.v. \int \frac{dy}{x-y}\rho_t(y) \right)\rho_t(x)\right), \label{eq:McKV0} \EEQ
in a weak (i.e. distribution) sense, where $p.v. \int \frac{dy}{x-y}\rho_t(y)$ is a principal value integral. 
\end{itemize}
The equilibrium measure $\rho_{eq}$, defined as the solution of the integral equation
(traditionally called: {\em cut equation})
\BEQ \frac{\beta}{2} p.v. \int \frac{dy}{x-y} \rho_{eq}(y)=V'(x) \label{eq:cut-eq1},
\EEQ  cancels the right-hand side 
of (\ref{eq:McKV0}), as is readily checked.

\medskip\noindent A complex Burgers-like PDE for the Stieltjes transform of $X_t$ 
\BEQ U_t(z):=\int \frac{1}{x-z} X_t(dx), \qquad z\in\C\setminus\R \EEQ
is easily derived \cite{RogShi,Isr} from (\ref{eq:McKV0}),
assuming $V$ to be polynomial,

\BEQ \frac{\partial U_t}{\partial t}=\frac{\partial}{\partial z} \left(
\frac{\beta}{4} (U_t(z))^2+ V'(z)U_t(z)+T_t(z) \right),  \label{eq:T}\EEQ
 where   
\BEQ T_t(z):=\int \frac{V'(x)-V'(z)}{x-z} \, X_t(dx).\EEQ

\bigskip
\noindent
In  our recent article \cite{Unt1}, in large part based on a previous paper by Israelsson \cite{Isr} which dealt with the specific example of a harmonic potential, we introduced
a process $Y=(Y_t)_{t\ge 0}$ interpreted as {\em asymptotic fluctuation process}.  Let $Y^N_t:=N(X^N_t-X_t)$ be the rescaled fluctuation process
for finite $N$. Then it was proved that  $Y^N_t \overset{law}{\to} Y_t$ when $N\to\infty$,
where $(Y_t)_{t\ge 0}$ is the solution of a martingale problem, as can be briefly seen
as follows. First, It\^o's formula implies that
\BEQ d\langle Y_t^N,f_t\rangle=\half(1-\frac{\beta}{2}) \langle X_t^N,f''_t\rangle dt+
\frac{1}{\sqrt{N}} \sum_{i=1}^N f'_t(\lambda_t^i) dW_t^i \label{eq:Ito} \EEQ
if the test functions $(f_t)_{0\le t\le T}$, $f_t:\R\to\R$ solve the following linear PDE
  
\BEQ \frac{\partial f_t}{\partial t}(x)=V'(x)f'_t(x)-\frac{\beta}{4} \int \frac{f'_t(x)-f'_t(y)}{x-y} (X_t^N(dy)+X_t(dy)) \label{eq:intro-PDE-f} 
\EEQ
Substituting formally to $X^N$ its deterministic limit $X$ in the
r.-h.s. of (\ref{eq:intro-PDE-f}), one gets an equation which is
the asymptotic limit of (\ref{eq:intro-PDE-f}) in the limit $N\to\infty$, namely,
\BEQ \frac{\partial f_t}{\partial t}(x)=V'(x)f'_t(x)-\frac{\beta}{2} \int \frac{f'_t(x)-f'_t(y)}{x-y} X_t(dy) \label{eq:intro-PDE-f-asymptotic} 
\EEQ
The main task in \cite{Unt1} consists in proving that eq. (\ref{eq:intro-PDE-f},
\ref{eq:intro-PDE-f-asymptotic}) is akin to a transport equation on the cut complex
plane $\C\setminus\R$. In the harmonic case (i.e. when $V$ is quadratic), then 
the solution of, say, (\ref{eq:intro-PDE-f-asymptotic}) at time $t$ with terminal condition
$f_T(x)\equiv \frac{c}{x-z},\ c\in\C,z\in\C\setminus\R$ is equal to $\frac{c_t}{x-z_t}$
where $z_t\equiv a_t+\II b_t$ and 
\BEQ \frac{da_t}{dt}=\frac{\beta}{2} \Re U_{t}(z_t) +V'(a_t), \qquad \
 \frac{db_t}{dt}=\frac{\beta}{2} \Im U_t(z_t)+ V''(a_t)b_t  \nonumber \EEQ 
\BEQ  \frac{dc_{t}}{dt}= \Big[\frac{\beta}{2}   U'_t(z_t) +2 V''(a_t) 
\Big] c_{t}.    \label{eq:char0} \EEQ
Thus the solution of (\ref{eq:intro-PDE-f-asymptotic}) may be represented formally as 
\BEQ \int da\, db\, h_t(a,b) \frac{1}{x-z}, \label{eq:density}
\EEQ
 where $h_t(a,b):=c_t \del(a-a_t)\del(b-b_t)$, interpreted as a {\em density} on
 $\C\setminus\R$, is obtained by "pushing" $h_T(a,b):=\del(a-a_T)
\del(b-b_T)$ along the above characteristics, or equivalently, by solving the
associated 
transport equation generated by the time-dependent operator 
\BEQ {\cal L}(t):=(\frac{\beta}{2} \Re U_{t}(z) +V'(a))\partial_a+ (\frac{\beta}{2} \Im U_t(z)+ V''(a)b )\partial_b+\frac{\beta}{2}   U'_t(z) +2 V''(a). \label{intro:L} 
\EEQ

\medskip\noindent     Considering instead some
arbitrary terminal condition and potential $V$, a similar formula holds, where the time-evolution  is given up to a bounded perturbation by a transport operator  whose
characteristics are as (\ref{eq:char0}) plus some extra term depending on $V'''$.
Then (at least formally), It\^o's formula (see \cite{Isr}, p. 29) makes it possible to
find  the Markov kernel in the limit $N\to\infty$. Namely, if $f_t$ be the solution of (\ref{eq:intro-PDE-f-asymptotic}) with terminal condition $f_T$, and $\phi_{f_t}(Y_t):=e^{\II\langle Y_t,f_t\rangle}$,
  \BEQ \esper[\phi_{f_T}(Y_T)\big| {\cal F}_t]=\esper[\phi_{f_t}(Y_t)] \exp\left( \half
\int_t^T  \left[ \II (1-\frac{\beta}{2}) \langle X_s,f''_s\rangle -
\langle X_s,(f'_s)^2\rangle \right] ds  \right). \label{eq:intro-Gaussian0} \EEQ
Eq. (\ref{eq:intro-Gaussian0}) was proved for  general potentials in our previous article \cite{Unt1}.
Now, letting $f_T(x):=\sum_{k=1}^n \frac{c_T^k}{x-z_T^k}$, $z_T^k\in\C\setminus\R, k=1,\ldots,n$  vary in dense subspace of $L^1(\R)$, this martingale problem is solved in
Bender \cite{Ben} in the case of a harmonic potential using an explicit computation of the characteristics (\ref{eq:char0}). Such is the present state of the art.


\subsection{Main results}


We prove in this article two types of results. We shall
generally assume that 
$V$ is {\em polynomial} and {\em strictly convex}, though the reader will
also find weaker sets of  hypotheses, depending on the paragraph.

\bigskip\noindent
{\bf (A)} The first series of  results regards the {\em Mc Kean-Vlasov equation} (\ref{eq:McKV0}).  
Little is known about it in general; the arguments in Li-Li-Xie \cite{LiLiXie} (see in particular
Theorem 1.3) simply prove that it
admits a unique solution in $C([0,T],{\cal P}(\R))$, which is constructed as weak limit
of the sequence of stochastic processes $t\mapsto Y_N(t)$. Unicity is proved using
decrease of Wasserstein distance between two arbitrary solutions. A classical large-deviation argument (reviewed here) implies under our hypotheses a bound on the support of the measure $\rho_t$;
in particular, $\rho_t$ is compactly supported. 

\medskip\noindent Our first result is a {\em regularity result}: assuming that the
analytic function $z\mapsto U_0(z), z\in\Pi_+:=\{\Im z>0\}$ extends to a continuous
function on the closure $\Pi_+\cup\R$ of the upper half-plane, we prove
that the same property holds for $U_t, t\ge 0$; see Theorem \ref{th:regularity}. Hence in particular (by Plemelj's formula),
the density $\rho_t(\cdot)=\frac{1}{2\II\pi}(U_t(\cdot+\II 0)-U_t(\cdot-\II 0))$ is a continuous function for every $t\ge 0$.

\medskip\noindent Our second result concerns the  {\em support}. We explain how to obtain the "external support"
$[a_t,b_t]$ of $\rho_t$, i.e. the intersection of all intervals $[a,b]$  such that 
$\langle \rho_t,\phi\rangle\equiv 0$ for every test function $\phi$ with compact support
$\subset\R\setminus[a,b]$.  (This implies that $\supp(\rho_t)\subset[a_t,b_t]$ but
not the reverse inclusion $[a_t,b_t]\subset\supp(\rho_t)$.)  The external support is characterized, see eq. (\ref{eq:b0*}) and (\ref{eq:atbt}),  in terms of
characteristics of the generalized complex Burgers equation (\ref{eq:T}) -- not
surprisingly closely related to (\ref{eq:char0}) --  which are half-explicit in general and can be obtained in  closed form in
 various cases, including for equilibrium dynamics or  when $V$ is harmonic. On the other hand, we do not prove any formula for the
 support itself. In particular, though under our hypotheses (more specifically, because
 $V$ is convex) the support of the equilibrium density  is a connected interval, we
 cannot exclude, even if $\supp(\rho_0)$ is connected, that  e.g. 
 $\supp(\rho_t)=[a,c]\cup[d,b]$ with $a<c< d<b$
for some  $t>0$.

\medskip\noindent Much stronger results have been
proved by P. Biane \cite{Bia} in the Hermite case, namely,
for $\beta=2$ and, say, 
 $V(x)=\frac{x^2}{2}$ (harmonic potential), see \S \ref{subsection:reg} {\bf A.} for more
 precise statements concerning regularity. Also, the number
 of connected components is shown to be decreasing with time, so that  the above hypothetical behavior can be excluded.

\bigskip\noindent
{\bf (B)}  The second series of results regards the {\em fluctuation process} $(Y_t)_{t\ge 0}$.  While the above characteristic equations can be solved explicitly only when $V$ is
harmonic (see Bender \cite{Ben}), yielding the covariance of the Stieltjes
transform  $({\mathcal{S}} Y_t)(z):=\langle Y_t,\frac{1}{\cdot-z}\rangle$ of the fluctuation process,
\BEQ \Lambda(t_1,z_1;t_2,z_2):=\Cov(({\mathcal{S}}Y_{t_1})(z_1),
({\mathcal{S}}Y_{t_2})(z_2)),
\EEQ
 their "trace" on the
boundary of the upper (or lower) half-planes can be solved for arbitrary $V$.
Then the covariance kernel $\Cov( Y_{t_1}(x_1),Y_{t_2}(x_2))$ is found by taking
boundary values $Y_{t_i}(x_i)=\frac{1}{2\II\pi}\Big ( ({\mathcal{S}}Y_{t_i})(x_i+\II 0)-({\mathcal{S}}Y_{t_i})(x_i-\II 0) \Big)$, $i=1,2$. Our most general result in this direction
is Theorem \ref{th:g-PDE}. A more explicit formula relying on Theorem \ref{th:g-PDE} is
Theorem \ref{th:g-stat} or Corollary \ref{th:g-stat} for equilibrium dynamics, see (\ref{eq:g-PDE-stat-quartic}) for  the specific case of a quartic
(Landau-Ginzburg type) potential.  The reader should compare the above results to
those obtained by M. Duits \cite{Dui} in a stochastic setting for fluctuations
of non-colliding processes, and by  N. Allegra, P. Calabrese, J. Dubail, J.-M. 
St\'ephan and J. Viti
\cite{Dub1},\cite{Dub2} in a condensed-matter context for the (real-time) propagator of the density field $\langle \rho(t_1,x_1)\rho(t_2,x_2)\rangle\equiv \langle (\psi^{\dagger}\psi)(t_1,x_1) (\psi^{\dagger}\psi)(t_2,x_2)\rangle$ of a one-dimensional Fermi gas submitted to a confining potential $V$.  Despite the difference of language, and the fact that an analytic continuation in time
is necessary to go from  one situation  to the other, both series of works come to
a similar conclusion.   Focusing on the quantum setting, and considering the low-lying spectrum of the underlying $N$-particle quantum Hamiltonian,  the authors   predict (and confirm by some numerical simulations) that (assuming the theory to be   free, i.e. Gaussian at large scale)   the time-evolution equation obtained  for the Wigner function in the semi-classical limit    is essentially correct in the large $N$ limit. The time-evolution equation for the chiral part of  the two-point
function  is then the same as ours (compare e.g. our
 equation (\ref{eq:short-distance-rho}) to eq. (6) in \cite{Dub2}), taking as input the equilibrium
density $\rho_{eq}$ computed by local-density approximation, see e.g.
discussion in section {\bf A.} of \cite{CamVic} or articles cited above. Then, in both situations, the
fluctuation/density field is  interpreted as a 2d Gaussian free field in a curved space
with metric tensor $ds^2= e^{2\sigma} dz\, d\bar{z}$, with  coordinate
transform $z=z(x,y)$ and conformal weight
$\sigma=\sigma(x,y)$ chosen by requiring that  $e^{\sigma(x,y)}dz=dx+\II \pi\rho_{eq}(x) dy$, which yields (\cite{Dub2}, eq. (20)): $z(x,y)=\frac{1}{\pi} ( G(x)+\II \pi y)$,
where $G(x):=\int \frac{dx}{\rho_{eq}(x)}$, in exact correspondence with  our Theorem \ref{th:g-stat}. Therefore its law may be obtained from
that of flat 2d Gaussian free field through a conformal transformation. The 
connection of our results to those  is however lost at that point,  since the single-time covariance
kernel Cov$(Y_{t}(x_1),Y_t(x_2))$ is (up to a simple scaling) independent of the potential, hence of $\rho_{eq}$. It would be interesting to obtain a deeper understanding of this difference.


\section{The Mc Kean-Vlasov equation}


We study in this section eq. (\ref{eq:McKV0}) indirectly through the time-evolution
of its Stieltjes transform 
\BEQ U_t(z):=\int dx\, \frac{\rho_t(x)}{x-z}, \qquad z\in\C\setminus\R. \EEQ
As shown in \cite{RogShi},\cite{Isr}, $U_t$ satisfies following generalized complex Burgers equation,
\BEQ \frac{\partial U_t(z)}{\partial t}=\frac{\partial}{\partial z}\big( \frac{\beta}{4} U_t^2(z)+ V'(z) U_t(z)+ T_t(z)\big),
\label{eq:U} \EEQ
where
\BEQ T_t(z):= \int dx\, \rho_t(x) \frac{V'(x)-V'(z)}{x-z}.\EEQ
When $V$ is {\em harmonic}, $T_t(z)$ is a constant, whence $T'_t(z)=\frac{d}{dz}
T_t(z)\equiv 0$. But in general,  $T_t$ is an unknown time-dependent quantity for which
an independent equation should be provided. For $V$ polynomial, however, say, deg$(V)=:2n$, $T_t(z)$
is easily seen \cite{Joh} to be some explicit polynomial in $z$ of order $\le 2n-2$, with
coefficients in the linear span of the  $2n-2$ first
moments of the unknown density $\rho_t$, namely, $T_t(z)=\sum_{k=0}^{2n-2}
c_k \int x^k \rho_t(x)\, dx$ for some constants $c_k=c_k[V]$.  Looking at the asymptotic
expansion of $U_t$ at infinity, $-U_t(z)\sim \frac{1}{z}+\Big( \int x\rho_t(x)\, dx
\Big) \frac{1}{z^2}+ \Big( \int x^2\rho_t(x)\, dx
\Big) \frac{1}{z^3}+\cdots$,  $T_t(z)$ may also be defined  (up to an additive constant) as minus the  part polynomial in $z$   of 
$V'(z)U_t(z)$, so that $\frac{\partial U_t(z)}{\partial t}=O(1/z^2)$ when $z\to\infty$, in coherence
with the leading term of the expansion, $-U_t(z)\sim_{z\to\infty} 1/z$. Projecting
(\ref{eq:U}) onto the linear subspace $\oplus_{k\ge 0} \C z^{-k-1}$ yields an infinite
system of coupled ODEs for the moments $\left(\int x^k \rho_t(x)\, dx\right)_{k\ge 0}$,
which in principle can be solved numerically on short time-intervals.

\medskip\noindent {\em Notation.} We let $\Pi_+:=\{z\in\C| \Im z>0\}$ be the upper half-plane, and $\bar{\Pi}_+:=\Pi_+\cup\R$
its closure.

\bigskip\noindent We make in this section the following

\medskip\noindent {\bf Assumptions.}  {\em
\begin{itemize}
\item[(i)] $V$ is a {\em strictly convex} polynomial of order $2n$, i.e. $\inf_{x\in\R} V''(x)\ge \alpha>0$;
\item[(ii)] $U_0\big|_{\Pi_+}$ extends to a {\em continuous function} $U_0:\bar{\Pi}_+\to\C$. 
\end{itemize} }

\noindent Since (by Plemelj's formula, see \S 5.1), $ \lim_{\eps\to 0^+} \Im U_0(x+\II\eps)=\pi\rho_0(x)$,  Assumption (ii)
implies in particular that $\rho_0$ is a continuous function.


\subsection{An example: scaling solution in the Hermite case}


In this paragraph, we assume that $\beta=2$ and $V(x)=\frac{x^2}{2}$, and look for
simple solution of (\ref{eq:U}) other than the constant solution $\rho_{eq}$.
By reference
to the underlying equilibrium unitary ensemble, we call this case the {\em Hermite case}.

\medskip\noindent
{\em Explicit formulas.} 
 The equilibrium density
corresponds to the semi-circle law, $\rho_{eq}(x)\equiv\frac{1}{\pi}\sqrt{2-x^2}
\, {\bf 1}_{|x|<\sqrt{2}}$,
with support $[-\sqrt{2},\sqrt{2}]$ and  Stieltjes transform $U_{eq}(z)\equiv -z+\sqrt{z^2-2}$ continuously extending to the real line,
\BEQ U_{eq}(x\pm\II 0)=-x\pm \II \sqrt{2-x^2} \qquad (|x|<2), \qquad
U_{eq}(x\pm\II 0)=-x+\sqrt{x^2-2} \qquad (|x|>2).\EEQ
Taking the boundary value, $\frac{1}{2\II\pi}(U_{eq}(x+\II 0)-U_{eq}(x-\II 0))\equiv
\frac{1}{\pi} \Im U_{eq}(x+\II 0)$, yields
$\rho_{eq}(x)$; the functions $U_{eq}(x\pm \II 0)$ are real-valued on $\R\setminus[-\sqrt{2},
\sqrt{2}]$ and $U_{eq}(\bar{z})=\overline{U_{eq}(z)}$, hence (by Schwarz's extension lemma)
$U_{eq}$ extends to a holomorphic function (still called $U_{eq}$) on the cut plane $\C\setminus[-\sqrt{2},\sqrt{2}]$. 
Note that $U'_{eq}$ is singular in the neighbourhood of
the ends of the support, $\pm\sqrt{2}$; namely, $U'_{eq}(\pm(\sqrt{2}+\eps))\sim_{\eps\to 0^+} \pm
c/\sqrt{\eps}$ ($c>0$).

\medskip\noindent {\em Scaling solution.  
Assume} that $\rho_0(x):=\frac{1}{s}\rho_{eq}(x/s)$ $(s>0)$, or equivalently, $U_0(z):=
\frac{1}{s} U_{eq}(z/s)$. Then we use the following Ansatz,
\BEQ U_t(z)\equiv \frac{1}{s(t)}U_{eq}(z/s(t)) \EEQ
for some unknown scaling function $t\mapsto s(t)$, corresponding to a time-dependent support
$[-s(t)\sqrt{2},s(t)\sqrt{2}]$.  From (\ref{eq:U}), we obtain
for the stationary solution $(U_0(z)+z)U'_0+U_0=0$. Hence
\BEA 0 &=& \frac{dU_t(z)}{dt}-(U_t(z)+z) U'_t(z)-U_t(z) \nonumber\\
&=&  
\frac{1}{s^2(t)} \left\{ (-U_{eq}(\frac{z}{s(t)})+\frac{z}{s(t)} U'_{eq}(\frac{z}{s(t)})) \dot{s}(t) - (\frac{1}{s(t)}U_{eq}(\frac{z}{s(t)}) + z) U'_{eq}(\frac{z}{s(t)}) - s(t) U_{eq}(\frac{z}{s(t)}) \right\} \nonumber\\
&=& -\frac{1}{s^2(t)}(U_{eq}(\frac{z}{s(t)})+ \frac{z}{s(t)} U'_{eq}(\frac{z}{s(t)}))
 \left\{ \dot{s}+s-\frac{1}{s} \right\}.
 \EEA
Hence our Ansatz is correct provided we choose $s(t)$ to be the solution of the ode
$\dot{s}=\frac{1}{s}-s$, namely,
\BEQ s(t)\equiv \sqrt{1+e^{-2t} (s^2(0)-1)}.\EEQ
Equivalently, $\frac{s^2(t)-1}{s^2(0)-1}=e^{-2t}$, which means that the "radius"
$b_t:=\sqrt{2}\, s(t)$ converges exponentially fast and monotonously to its equilibrium value, 
$\sqrt{2}$. 


\subsection{Regularity}  \label{subsection:reg}


As proved in our previous article \cite{Unt1} -- extending uniform-in-time moment
bounds proved in \cite{And} in the harmonic case --, there exists $R=R(T)$ and
$c,C>0$  such
that, for all $N\ge 1$, $\proba[\sup_{0\le t\le T} \sup_{i=1,\ldots,N} |\lambda_t^{N,i}|
>R]\le Ce^{-cN}$ (see Proposition \ref{prop:LD}). Using Borel-Cantelli's lemma, one immediately deduces the following: for any test function $f:\R\to\R$ with support $\subset
B(0,R)^c$,  $\langle \rho_t,f\rangle=\lim_{N\to\infty} \langle X^n_t,f\rangle=0$ a.s.
Thus $\supp(\rho_t)\subset[-R,R]$ for every $t\le T$. In particular, for every $n=0,1,\ldots$,  the function $t\mapsto \int x^n \rho_t(x)\, dx$ $(0\le t\le T)$  is bounded and
continuous; which implies in turn that $t\mapsto T'_t(z)$ is a polynomial in $z$ depending
continuously on $t$. 

\bigskip
\noindent Our main result in this subsection is

\begin{Theorem}  \label{th:regularity}
Under the Assumptions of section 2, $U_t\Big|_{\Pi_+}$ extends to a contiuous function
on $\bar{\Pi}_+$ for every $t\ge 0$. In particular, $x\mapsto \rho_t(x)$ is a continuous
function for every $t\ge 0$.
\end{Theorem}

\noindent\bigskip {\bf A. (Case of a harmonic potential).}

\medskip\noindent Then  $\frac{d}{dz} T_t(z)\equiv 0$ and so (\ref{eq:U}) is a closed
equation for $U_t$ which can be solved on $\C\setminus\R$, where it is analytic, using the method of characteristics. We shall use this to derive the evolution of the support.

\medskip\noindent {\em Characteristics.} For definiteness we choose $V(x)=\frac{x^2}{2}$. Let $Z_t(z_0)$ be the solution at time $t\ge 0$ of the following differential
equation,
\BEQ  \frac{dz}{dt}=-\frac{\beta}{2} U(t,z(t))-z(t), \qquad z(0)=z_0\in\Pi^+.  \label{eq:edo}
\EEQ
Letting $C(t):=-U(t,z(t))$ and substituting into (\ref{eq:U}) yields $\frac{d}{dt} C(t)=C(t)$, solved as $C(t)=e^t C(0)$.  Differentiating (\ref{eq:edo}) yields 
\BEQ \ddot{z}=\frac{\beta}{2} \dot{C} - \dot{z}=z \EEQ
with $z(0)=z_0, \dot{z}(0)=-\frac{\beta}{2}U_0(z_0)-z_0$
hence
\BEQ Z_t(z_0)=z_0 \ch t - \big[ \frac{\beta}{2} U_0(z_0)+z_0\big] \sh t=
z_0 e^{-t} -\frac{\beta}{2} U_0(z_0) \sh t. \label{eq:Z} \EEQ
Since $\Im U_0(z_0)\ge 0$ by (\ref{eq:positivity}) for $z_0\in\Pi_+$, $t\mapsto \Im Z_t(z_0)$ decreases and the
characteristics may eventually cross the real axis, after which the characteristic method
makes  no sense because of the discontinuity. So we decide to {\em kill} 
characteristics as sooon as they cross the real axis.

\noindent Let $t_{max}(z_0):=\inf\{t>0 \ |\ Z_t(z_0)\in\R\}\in (0,+\infty]$; for every
$T<t_{max}(z_0)$, there exists a neighbourhood ${\cal B}(z_0)$ of $z_0$ in $\Pi_+$ that is mapped
inside $\Pi_+$. Hence characteristics (\ref{eq:edo}) started from
${\cal B}(z_0)$ are well-defined up to time $T$, and define for every $t\le T$ a one-to-one mapping into
a time-dependent region $Z_t({\cal B}(z_0))\subset\Pi_+$. 
{\em Denote by $\phi_t:Z_t({\cal B}(z_0))\to {\cal B}(z_0)$ the inverse mapping}, $\phi_t(z):=Z_t^{-1}(z)$. Then
\BEQ U_t(z)=e^t U_0(\phi_t(z)), z\in Z_t({\cal B}(z_0)). \EEQ
 Solving instead {\em backwards}
in time, one gets
\BEQ \phi_t(z)=ze^t+\frac{\beta}{2} U_t(z)\sh t. \label{eq:phit} \EEQ
Since $\Im U_t(z)\ge 0$, it is apparent from (\ref{eq:phit}) that $\phi_t:\Pi_+\to
\Pi_+$, with $\Im \phi_t(z)\ge \Im z$; this can be deduced, even without knowing the
explicit formula (\ref{eq:phit}), from (\ref{eq:edo}), since $-\frac{dz}{dt}\in\Pi_+$
as long as $z(t)\in\Pi_+$. Let 
\BEQ \Pi_t:=\phi_t(\Pi_+) \EEQ
and $\bar{\Pi}_t\subset\bar{\Pi}_+$ its closure in $\bar{\Pi}_+$.  Since
$\Pi_t=\{z\in\Pi_+\ Z_s(z)\in\Pi_+, 0\le s\le t\}$, it is clear that the family of
regions $(\Pi_t)_{t\ge 0}$ is decreasing for inclusion, i.e. $\Pi_t\supset\Pi_T$ for
$T\ge t$.  If
$w_n:=\phi_t(z_n)$, $z_n\in\Pi_+$ is a sequence in $\Pi_t$ converging to $w$,
then $z_n=Z_t(w_n)\to Z_t(w)$ by (\ref{eq:Z}) since $U_0$ is continuous on $\bar{\Pi}_+$. Furthermore, if $|w_n|\to\infty$, then $|z_n|\sim e^{-t} |w_n|\to\infty$. Thus (see
Rudin \cite{Rud}, Theorem 14.19) the map $\phi_t$ extends to a homeomorphism
$\bar{\Pi}_+\to \bar{\Pi}_t$, while the boundary $\partial \bar{\Pi}_t$ is a
Jordan curve. {\em Hence $U_t:z\mapsto e^t U_0(\phi_t(z))$
extends to a continuous function on $\bar{\Pi}_+$.}

\bigskip\noindent Specifically when $\beta=2$ {\em (Hermite case)}, the density at time $t$ may be interpreted as
the free convolution of the time-zero density by a semi-circular law. P. Biane \cite{Bia} proves then much more.
First, whatever the initial condition, the measure at time $t>0$ has a continuous density. Then, the density $\rho_t$ is proved to
be analytic on the open subset $\{\rho_t>0\}:=\{x\in\R\ |\ 
\rho_t(x)>0\}$. Furthermore,
$\overline{\{\rho_t>0\}}$ is the support of the measure at time $t$, and  $\rho_t(x)=O((d(x,\{\rho_t>0\}^c)^{1/3})$ for
$x\in\{\rho_t>0\}$, where
$d(x,\{\rho_t>0\}^c)$ is the distance to the complentary
set. This a priori surprising $1/3$-H\"older exponent gives
the correct behavior of $\rho_t$ at a point $x$ where two components
of the support merge at time $t$.

\bigskip\noindent {\em A simple example.} Assume $U_t=U_{eq}$, $t\ge 0$. Then
$\phi_t(z)=ze^t+(-z+\sqrt{z^2-2})\sh t= z\ch t +\sqrt{z^2-2}\, \sh t$, whence
(for $|x|<2$)
$\phi_t(x+\II 0)\equiv a(x)+\II b(x)$, with $a(x)=x\ch t,b(x)=\sqrt{2-x^2}\, \sh t$.
Thus $\partial \bar{\Pi}_t$ is the union of $(-\infty,-\ch t]\cup[\ch t,+\infty)$ with the semi-ellipse defined by the equation $\{\frac{a^2}{\ch^2 t} + \frac{b^2}{\sh^2 t}=2,\ b\ge 0\}$. This makes it plain enough that (somewhat counter-intuitively) characteristics
 {\em do not} follow the time-evolution of the support or the singularities of
 $U_t$ on the real axis (see next subsection for more). 


\bigskip\noindent {\bf B.  General case}

\medskip\noindent The general case is similar, except that the time evolution of 
the  
$(2n-2)$ first moments of the density must be determined independently.  Namely,  instead of 
(\ref{eq:edo}), we consider the generalized characteristics $Z_t(z_0)$, solution of the
o.d.e.
\BEQ \frac{dz}{dt}=-\frac{\beta}{2} U(t,z(t))-V'(z(t)), \qquad z(0)=z_0.  \label{eq:edo-bis}
\EEQ
and the mapping $\phi_t\equiv Z_t^{-1}$.
Letting $C(t):=-U(t,z(t))$ and substituting into (\ref{eq:U}) yields $\frac{d}{dt} C(t)=V''(z(t))C(t)-T'_t(z(t))$, solved as 
\BEQ C(t)=A_0^t C(0)-\int_0^t dt'\, A_{t'}^t T'_{t'}(z(t')), \qquad
A_{t'}^t:=\exp\left( \int_{t'}^t ds\, V''(z(s))\right).  \label{eq:C(t)-bis} \EEQ
  Differentiating (\ref{eq:edo-bis}) yields 
\BEQ \ddot{z}=V''(z) V'(z)-\frac{\beta}{2} T'_t(z)  \label{eq:ddotz} \EEQ
with initial condition
\BEQ z(0)=z_0, \qquad \dot{z}(0)=-\frac{\beta}{2}U_0(z_0)-V'(z_0) \EEQ
whence 
\BEQ \dot{z}:=\pm\sqrt{ (V'(z))^2-\beta (T_t(z)-T_0(z_0))+\beta U_0(z_0)(\frac{\beta}{4} U_0(z_0)+V'(z_0))}.  \label{eq:z-sqrt} \EEQ
 Solving for $T_t$ by some  independent means (e.g. numerically), (\ref{eq:z-sqrt}) can
be solved numerically for short time  knowing $U_0$ (and even by quadrature when $T_t$
is constant, e.g. for equilibrium dynamics). However (due to the
multi-valuedness of the square-root function on $\C$), eq. (\ref{eq:z-sqrt})
stops making sense in general when the function inside the square-root vanishes. On the
other hand, an
unambiguous definition may be given in terms of the second-order differential equation (\ref{eq:ddotz}), in its matrix form 
\BEQ \frac{d}{dt} \left(\begin{array}{c} z\\ \dot{z}\end{array}\right)=\left(\begin{array}{c} \dot{z}\\ V''(z) V'(z)-\frac{\beta}{2} T'_t(z)  \end{array}\right).
\label{eq:ddotz-vec} \EEQ 

 Writing $V'(z)\sim_{z\to\infty} c_n z^{2n-1}+\ldots$, we get for $0<b<1$: $\Im V'(a+\II b)\sim_{a\to\infty} (2n-1)c_n a^{2n-2}b$, whence there exists
$a_{max}\ge 0$ such that:
\BEQ \Big(0<b<1,|a|\ge a_{max}\Big)\Rightarrow  \Im V'(a+\II b)>0.\EEQ 
On the other hand, since $V$ is strictly convex, there exists $b_{max}\in (0,1)$ such that
\BEQ \Big(0<b<b_{max},|a|\le a_{max}\Big)\Rightarrow
\Re V''(a+\II b)>0; \EEQ
for such $a,b$ one thus gets $V'(a+\II b)-V'(a)=\II \int_0^b V''(a+\II y)\, dy\in \Pi_+$.
 Thus (see (\ref{eq:edo-bis})) $-\frac{dz}{dt}\in\Pi_+$ as in the harmonic
case, providing one restricts to the strip $\Im z\in(0,b_{max})$.  The rest of the
argument proceeds as in the previous subsection if one restricts to characteristics 
included either in  $[-a_{max},a_{max}]\times[0,b_{max}]$ or in $(\R\setminus
[-a_{max},a_{max}])\times [0,1]$.  {\em Hence, letting $z_0\equiv \phi_t(z)$ so that
$z(t')=Z_{t'}(\phi_t(z))=\phi_{t-t'}(z)$, 
\BEQ U_t:z\mapsto A_0^t U_0(\phi_t(z))+\int_0^t dt'\, A_{t'}^T T^t_{t'}(\phi_{t-t'}(z)),
\label{eq:UAAT}
\EEQ
 see
(\ref{eq:edo-bis},\ref{eq:C(t)-bis}) with $z(0)=\phi_t(z)$,
extends to a continuous function on $\bar{\Pi}_+$,} proving Theorem \ref{th:regularity} in
whole generality.


\subsection{Support} \label{subsection:support}


  In
this paragraph we study the time evolution of the {\em external support} $[a_t,b_t]$
defined as the intersection of all intervals $[a,b]$  such that 
$\langle \rho_t,\phi\rangle\equiv 0$ for every test function $\phi$ with compact support
$\subset\R\setminus[a,b]$.
 Using the characteristics introduced in the previous subsection, we shall be able
to give a defining formula for $a_t, b_t$ $(t\ge 0)$.

\medskip\noindent
Exactly as in the example developed in \S 2.1, and for the same reasons, the function
$U_0$ has a maximal analytic extension  to the cut plane $\C\setminus[a_0,b_0]$, which
is real-valued and real-analytic on $\R\setminus[a_0,b_0]$. Thus the characteristics  $
t\mapsto Z_t(x_0)$ issued from $x_0>b_0$, as defined by (\ref{eq:edo-bis}), is well-defined
and real-valued for $t$ small enough. As long as the characteristics $(z_s)_{0\le s\le t}$, 
$z_s:=Z_s(x_0)$ remains $\gg b_0$, i.e. for $x_0$ large enough, the dominant term inside the square-root 
in (\ref{eq:z-sqrt}) is $(V'(z))^2\sim (c_n z^{2n-1})^2$ ($c_n>0$),  the sign is
unambiguously a minus sign, $\dot{z}\approx -V'(z)$, and characteristics may not cross:
for $t\le T$ fixed and $b_{max}>b_0$ large enough,
the mapping $[b_{max},+\infty)\to \R, x_0\mapsto Z_t(x_0)$ is an increasing, real-analytic diffeomorphism on its
image. On the other hand, taking the derivative of (\ref{eq:ddotz-vec}) with
respect to the initial condition $\left(\begin{array}{c} x_0\\ \dot{x}_0\end{array}\right)$,
$\dot{x}_0=-\frac{\beta}{2} U(x_0)-V'(x_0)$, one can in
general only write
$Z'_t(x_0)$ in terms of some time-ordered exponential of matrices,
\BEQ  Z'_t(x_0)= \Big[  \overrightarrow{\exp} \Big( \int_0^t  ds\, \left( \begin{array}{cc}
0 & 1 \\ (V''V'-\frac{\beta}{2} T'_s)'(x_s)  & 0 \end{array} \right)  \Big)  \, \cdot\,  \left(\begin{array}{c}
1 \\ -\frac{\beta}{2} U'_0(x_0)-V''(x_0) \end{array}\right)\Big]_1,  \label{eq:Z't}
\EEQ
$[\, \cdot\, ]_1$=1st component, a complicated formula from which no general rule
to guess the possible vanishing of $Z'_t(x_0)$ can be expected.  Let us illustrate this on the simple Hermite case where $\beta=2$ and $V(x)=\frac{x^2}{2}$, and characteristics are explicit (see
{\bf A.} of last subsection).
When $x_0\to \infty$,  $Z_t(x_0) \sim e^{-t} x_0 +  \frac{\beta}{2} \sh t\, x_0^{-1}
+ O(x_0^{-2})$, hence in particular $x_0\mapsto Z_t(x_0)$ is increasing for 
$x_0$ large. On the other hand, one may expect that  $\frac{dZ_t(x_0)}{dx_0}=\ch t - \left[
\frac{\beta}{2} U'_0(x_0)+1\right] \sh t \to_{x_0\to b_0^+} -\infty$  for all $t>0$, which
does happen
e.g. when $U_0(z)=\frac{1}{s}U_{eq}(z/s)$ is a rescaling of the equilibrium solution
$U_{eq}$.

\medskip\noindent Define:
\BEQ b_0^*(t):=\sup\{x_0> b_0 \ |\  \min_{s\in[0,t]} Z'_s(x_0) \le 0\}. \EEQ
 By construction, $(b_0^*(t),+\infty)$ is the
largest interval of the form $(x_0,+\infty)$, $x_0\ge b_0$, such that
$Z_s:(b^*_0(t),+\infty)\to (Z_s(b^*_0(t)),+\infty)$ is a diffeomorphism for all $0\le s\le t$.  Then $Z_t\Big|_{(b^*_0(t),+\infty)}$ extends analytically on some complex neighbourhood ${\cal B}(b^*_0(t),+\infty)$ of $(b^*_0(t),+\infty)$ to
a conformal mapping with inverse $\phi_t$. Thus the function $U_t$ defined on
the image $Z_t({\cal B}(b^*_0(t),+\infty))$ by (\ref{eq:UAAT}) is a holomorphic solution
of (\ref{eq:U}). Hence $\supp(\rho_t)\subset(-\infty,b_t]$, where
\BEQ b_t:=Z_t(b^*_0(t)).\EEQ   
Also, $\min_{s\in[0,t]} Z'_s(b^*_0(t))=0$, so let $s_0:=\min\{s\in[0,t]\, |\, 
Z'_s(b^*_0(t))=0\}$. Then $Z'_s(b_0^*(t))>0$ for all $s<s_0$, so necessarily
$\frac{d}{dx_0} \dot{Z}_{s_0}(b_0^*(t))=\frac{d}{ds_0}Z'_{s_0}(b_0^*(t))\le 0$.  Since 
\BEA  &&\frac{d}{ds_0}Z'_{s_0}(x_0)=
\Big[  \left(\begin{array}{cc} 0 & 1 \\ 
(V''V'-\frac{\beta}{2} T'_{s_0})'(x_{s_0})  & 0 \end{array} \right) \ \cdot \nonumber\\
&&\qquad \cdot\  \overrightarrow{\exp} \Big( \int_0^{s_0}  ds\, \left( \begin{array}{cc}
0 & 1 \\ (V''V'-\frac{\beta}{2} T'_s)'(x_s)  & 0 \end{array} \right)  \Big)  \, \cdot\,  \left(\begin{array}{c}
1 \\ -\frac{\beta}{2} U'_0(x_0)-V''(x_0) \end{array}\right)\Big]_1 \nonumber\\
&&\ \  = \Big[  \overrightarrow{\exp} \Big( \int_0^{s_0}  ds\, \left( \begin{array}{cc}
0 & 1 \\ (V''V'-\frac{\beta}{2} T'_s)'(x_s)  & 0 \end{array} \right)  \Big)  \, \cdot\,  \left(\begin{array}{c}
1 \\ -\frac{\beta}{2} U'_0(x_0)-V''(x_0) \end{array}\right)\Big]_2, \nonumber\\
\EEA
the simultaneous  vanishing  of $Z'_{s_0}(b^*_0(t))$ and
$\frac{d}{ds_0}Z'_{s_0}(b_0^*(t))$ implies that the
two-component vector between square brackets $\Big[\ \cdot\
\Big]$ in (\ref{eq:Z't}) vanishes, which is impossible since
$\overrightarrow{\exp}\Big(\ \cdot\ \Big)$ is invertible. So, actually,
$\frac{d}{ds_0}Z'_{s_0}(b_0^*(t))<0$, which is contradictory with the definition of
$b_0^*(t)$ if $s_0<t$.  Hence $b^*_0(t)$ is also defined more simply as

\BEA b_0^*(t)&:=&\sup\{x_0> b_0 \ |\  Z'_t(x_0) \le 0\} \nonumber\\
&=& \inf\{x_0>b_0\ |\ Z_t(\cdot):(x_0,+\infty)\to(Z_t(x_0),+\infty)\ {\mbox{is a 
diffeomorphism}} \}. \nonumber\\ \label{eq:b0*}
\EEA 
Conversely, suppose $U_t$ were analytic at $b^*_0(t)$, then (\ref{eq:edo-bis}) would
imply that $Z'_t(b_0^*(t))\not=0$, a contradiction. 
Hence $\supp(\rho_t)\not\subset(-\infty,b_t-\eps)$ for any $\eps>0$.

\medskip\noindent {\em We now claim that the function $t\mapsto b_t$ is c\`adl\`ag, i.e.
right-continuous with left limits. Furthermore,  it doesn't have any positive jumps, 
i.e. $b_t\le\lim_{t'\to t,t'<t} b_{t'}$}. (On the other hand, we cannot exclude
negative jumps, with $\rho_{t'}\big|_{[b_t,b_{t^-}]}\to_{t'\to t,t'<t} 0$ pointwise). Namely, (i) $b_t\le\liminf_{t'\to t} b_{t'}$; otherwise
(by absurd), letting $b\in(\liminf_{t'\to t} b_{t'},b_t)$,  we would have $\int_{-\infty}^{b} \rho_{t_n}(x) dx\to 1>\int_{-\infty}^b
\rho_t(x)dx$, where $(t_n)_{n\ge 1}$ is a sequence such that $t_n\to t$, $b_{t_n}\to \liminf_{t'\to t}b_{t'}$ and 
$b_{t_n}<b$, which is incompatible  with the fact that
the measure $\rho_s(x) dx$ depends continuously on $s$; (ii) $\limsup_{t'\to t, t'>t}
b_{t'}\le b_t$, as follows from the characteristic method developed above; (iii) imagine (by absurd) that $b^-_{min}:=\liminf_{t'\to t,
t'<t} b_{t'}< b^-_{max}:=\limsup_{t'\to t,t'<t} b_{t'}$. Choose $b^-_{min}<b<b'<b^-_{max}$. Let $t_n\to t$ (resp. $t'_n\to t$) a sequence such
that $t_n,t'_n<t$ and $b_{t_n}\le b$ (resp. $b_{t'_n}\ge b'$). Then there exist characteristics moving by an amount $b'-b$ in arbitrary small time, which is contradictory
with  previous arguments.

\medskip\noindent We similarly define $a_0^*(t)$ by requiring that $(-\infty,a_0^*(t)]$ be the
largest interval of the form $(-\infty,x_0)$, $x_0\le a_0$, such that
$Z_t(\cdot):(-\infty,x_0)\to (-\infty,Z_t(x_0))$ is a diffeomorphism.
Then it follows from the above that
\BEQ [a_t,b_t]:=[Z_t(a_0^*(t)),Z_t(b_0^*(t))]  \label{eq:atbt} \EEQ
is the external support of $\rho_t$.

\bigskip\noindent 
{\em Let us illustrate this with the  example of the scaling solution of \S 2.1.} 
We find from (\ref{eq:Z})
\BEQ Z_t(z_0)=z_0 \ch t + \left[ \frac{2}{b_0^2} (z_0-\sqrt{z_0^2-b_0^2})-z_0\right]
\sh t,\qquad \frac{dZ_t(z_0)}{dz_0}= \ch t + \left[ \frac{2}{b_0^2} (1-
\frac{z_0}{\sqrt{z_0^2-b_0^2}})-1\right] \sh t  
\EEQ
The Jacobian $x_0\mapsto \frac{dZ_t(x_0)}{dx_0}$ vanishes for a single value 
$x_0^*(t)>b_0$, determined by
\BEQ \frac{x_0^*(t)}{\sqrt{(x_0^*(t))^2-b_0^2}}=1+\frac{b_0^2}{2}(\coth t-1),\qquad
\sqrt{(x_0^*(t))^2-b_0^2}=\left[ \frac{b_0^2}{4}(\coth t-1)^2+ (\coth t-1)\right]^{-1/2}.\EEQ
Easy but tedious computations yield
\BEQ 2\sh^2 t/((x_0^*(t))^2-b_0^2)=1+ (\frac{b_0^2}{2}-1)e^{-2t} \EEQ
\BEA  Z_t(x_0^*(t))&=& \sqrt{(x_0^*(t))^2-b_0^2} \left\{ (1+\frac{b_0^2}{2}(
\coth t-1))\ch t + \left[ \coth t-1-(1+\frac{b_0^2}{2}(\coth t-1))\right] \sh t \right\} \nonumber\\
&=& \sqrt{\frac{(x_0^*(t))^2-b_0^2}{\sh^2 t}} \left\{ (2-r_0^2)\sh t\,  \ch t + 
\frac{b_0^2}{2} (\ch^2 t + \sh^2 t) -2\sh^2 t\right\} \nonumber\\
&=& \sqrt{\frac{(x_0^*(t))^2-b_0^2}{\sh^2 t}} \left\{1+ (\frac{b_0^2}{2}-1)e^{-2t}
\right\} \nonumber\\
&=& \sqrt{2}\, s(t)
\EEA
as expected.


\section{Kernel of the fluctuation process}


We give in this section formulas for the distribution-valued covariance kernel
\BEQ g_{1,2}(t_1,x_1;t_2,x_2):=\Cov\Big( Y_{t_1}(x_1) ,Y_{t_2}(x_2)\Big) \EEQ
 of the asymptotic
fluctuation process $(Y_t)_{t\ge 0}$. The proof is indirect. First we obtain an
evolution equation for the Stieltjes transformed covariance kernels
 \BEQ g^{\eps_1,\eps_2}_{1,2}(t_1,x_1;t_2,x_2):=\lim_{y\to 0^+} \Lambda(t_1,x_1+ \II\eps_1 y ;t_2,x_2+\II \eps_2 y), \qquad \eps_1,\eps_2=\pm 
 \label{eq:g-Lambda} \EEQ
which are the boundary values of the kernel $\Lambda:(\C\setminus\R)\times(\C\setminus\R)
\to \C$ defined by
\BEQ \Lambda(t_1,z_1;t_2,z_2):=\Cov\Big( ({\cal S}Y_{t_1})(z_1), ({\cal S}Y_{t_2})(z_2) \Big).  \label{eq:3.3}
\EEQ
where $\Cov(Z_1,Z_2)$ for two complex-valued random variables $Z_1=X_1+\II Y_1,Z_2=X_2+\II Y_2$ means $\Cov(X_1,X_2)-\Cov(Y_1,Y_2)+\II(\Cov(X_1,Y_2)+\Cov(X_2,Y_1))$.   In (\ref{eq:3.3}), ${\cal S}Y_t:\C\setminus\R\to\C$ is the Stieltjes transform of $Y_t$,
$({\cal S}Y_t)(z):=\langle Y_t,\frac{1}{\cdot-z}\rangle$. 
Then we use Plemelj's formula, $\frac{1}{x-(x_j+\II 0)}-\frac{1}{x-(x_j-\II 0)}=2\II \pi \del(x-x_j)$, $j=1,2$, and obtain
\BEA  g_{1,2}(t_1,x_1;t_2,x_2) &=&
-\frac{1}{4\pi^2} \Big[g^{+,+}_{1,2}-g^{-,+}_{1,2}-g^{+,-}_{1,2}+g^{-,-}_{1,2}\Big] (t_1,x_1;t_2,x_2)  \nonumber\\
&=& -\frac{1}{2\pi^2}\  \Re \Big[ g^{+,+}_{1,2}-g^{+,-}_{1,2}\Big](t_1,x_1;t_2,x_2)
 \label{eq:gggg}
  \EEA
All these formulas are to be understood in a distribution sense. 

\medskip\noindent Though we are not able to solve (\ref{eq:intro-Gaussian0}) for an arbitrary test
function $f_T$, it turns out that the limiting evolution equation for $f_T(x):=
\sum_{k=1}^n \frac{c_T^k}{x-z_T^k}$ when
$\Im z_T^k\to 0^+$ (see Introduction) is an explicit transport equation, which is the key to the PDE
we obtain for the kernel $g^{\pm,\pm}$; see Theorem \ref{th:g-PDE}. This PDE can be solved in terms of the
characteristics (see (\ref{eq:g-PDE-solution})). In the stationary
case one gets a more explicit formula (see Theorem \ref{th:g-stat} and
Corollary \ref{cor:g-stat}). 

\medskip\noindent We end this section with the interesting case of a quartic potential,
$V(x)=\frac{1}{4} t^4+\frac{c}{2} t^2+d$ ($c>0$), for which computations can be made
totally explicit (see eq. (\ref{eq:g-PDE-stat-quartic})).


\subsection{General framework}


We collect here those notations and results proved in our previous article \cite{Unt1} 
which are necessary for the present study.


\subsubsection{Assumptions}


Our Assumptions in this section are of three different types.

\medskip\noindent {\bf Assumptions on the potential.}

\medskip\noindent {\em We assume that $V$ is convex and $C^{11}$.}

\medskip\noindent  The convexity assumption on $V$ is essential for the convergence of the finite $N$-density
to the solution $\rho_t$ of the Mc Kean-Vlasov equation, see \cite{LiLiXie}, and for  Johansson's universal
formula for equilibrium fluctuations to apply \cite{Joh}, see \S 3.4 below. The extra
regularity assumptions on $V$ have been used in \cite{Unt1} for semi-group estimates 
and in some perturbation arguments. Later on (see end of \S 3.3, and \S 3.4),  we
shall further assume that $V$ extends analytically to an entire function 
$V:\C\to\C$ in order to get more explicit formulas.

\medskip\noindent {\bf Assumptions on the initial  measure.}

\medskip\noindent Let $\mu^N_0=\mu_0(\{\lambda^i_0\}_i)$ be the initial measure of the stochastic process
$\{\lambda^i_t\}_{t\ge 0, i=1,\ldots,N}$, and $X_0^N:=\frac{1}{N}\sum_{i=1}^N
\del_{\lambda^i_0}$ be the initial empirical measure. Since $N$ varies,
we find it useful here to add an extra upper index $(\lambda^{N,i}_0)_{
i=1,\ldots,N}$ to denote the initial condition of the process for a given value of $N$.   {\em We assume that:

\medskip
\begin{itemize}
 \item[(i)]  (large deviation estimate for the initial support) there exist some constants $C_0,c_0,R_0>0$  such that, for every $N\ge 1$,
\BEQ \proba[ \max_{i=1,\ldots,N} |\lambda_0^{N,i}| > R_0] \le 
C_0 e^{-c_0 N }. \label{eq:support-bound0} \EEQ
 \item[(ii)]  $X_0^N\overset{law}{\to}
\rho_0(x)\, dx$ when $N\to\infty$, where $\rho_0(x)$ is a deterministic measure;
\item[(iii)] (rate of convergence)  
\BEQ \left(\esper[|U_0^N(z)-U_0(z)|^2]\right)^{1/2}=O(\frac{1}{N b}) \label{hyp:initial-measure} \EEQ
 for $z=a+\II b\in\C\setminus\R$, where $U_0(z):=\int dx\, \frac{\rho_0(x)}{x-z}$ is the
Stieltjes transform of $\rho_0$.

\end{itemize}}

\medskip\noindent 
As proved in \cite{Unt1},  
the initial large deviation estimate (i) implies a uniform-in-time large deviation estimate for the
support of the random point measure:

\begin{Proposition} (see \cite{Unt1}, Lemma 5.1)  \label{prop:LD}
Assume (i) holds for some constants $R_0,c_0,C_0>0$.
 Let $T>0$. There exists some radius $R=R(T)$ and constant $c$, depending on $V$ and $R_0,c_0$ but uniform in $N$,  such that
\BEQ \proba\left[ \sup_{0\le t\le T} \sup_{i=1,\ldots,N} |\lambda_t^{N,i}| > R\right]
\le C e^{-cN}. 
\EEQ
\end{Proposition}

\medskip\noindent Finally, as in section 2, we add a 

\medskip\noindent {\bf Regularity assumption on the initial density.}

\medskip\noindent {\em We assume that the Stieltjes transform $U_0\Big|_{\Pi_+}$
 of the initial density $\rho_0$  on the upper-half plane extends to a continuous function
 $\bar{\Pi}_+\to\C$.}
 
\medskip\noindent Though this Assumption is probably unnecessary, it is natural, holds true in
all examples treated  below, and allows stating  convergence
results in a stronger sense.   


\subsubsection{Summary of results}


All results presented here come from our previous article \cite{Unt1}.

\medskip\noindent {\bf Notation.} Generally speaking and without further mention, if
$z,z_1,z_2,\ldots$ are complex numbers, then we write
$z=a+\II b,z_1=a_1+\II b_1,z_2=a_2+\II b_2,\ldots$ their decomposition into real/imaginary part.

\begin{Definition}[Sobolev spaces] \label{def:Sobolev}
Let $H_n:=\{f\in L^2(\R)\ |\ ||f||_{H_n}<\infty\}$ $(n\ge 0)$, where 
$||f||_{H_n}:=\left(\int d\xi\, (1+|\xi|^2)^{n} |{\cal F}f(\xi)|^2\right)^{1/2}$, and 
$H_{-n}:=(H_n)'$ its dual. 
\end{Definition}

\medskip
\noindent The measure-valued process 
\BEQ Y^N:=N(X_t^N-X_t)  \label{eq:finite-N-process} \EEQ
 has been shown in \cite{Unt1}  to converge in $C([0,T],H_{-14})$:

\begin{Proposition}[Gaussianity of limit fluctuation process] (see \cite{Unt1}, Main
Theorem)  \label{prop:Gaussian}
 Let $Y_t^N$ be the finite $N$ fluctuation process (\ref{eq:finite-N-process}). Then:
\begin{enumerate}
\item $Y^N \overset{law}{\to} Y$ when $N\to\infty$, where $Y$ is a  Gaussian process.
More precisely, $Y^N$ converges to $Y$ weakly in $C([0,T],H_{-14})$; 
\item  let $\phi_f(Y_t):=e^{\II \langle Y_t,f\rangle}$.
Then 
\BEQ \esper[\phi_{f_T}(Y_T)\big| {\cal F}_t]=\phi_{f_t}(Y_t) \ 
 \exp\left( \half
\int_t^T  \left[ \II (1-\frac{\beta}{2}) \langle X_s,(f_s)''\rangle -
\langle X_s,((f_s)')^2\rangle \right] ds  \right)  \label{eq:intro-Gaussian} \EEQ
where $(f_s)_{0\le s\le T}$ is the solution of the asymptotic equation (\ref{eq:intro-PDE-f-asymptotic}).
 
\end{enumerate} 
\end{Proposition}

\noindent The main point of the proof has been to rewrite the evolution equation for $(f_t)_{0\le t\le T}$ in terms of a "quasi"-transport operator on functions on the upper half-plane. Let us briefly recapitulate how this is done.

\begin{Definition}[Stieltjes transform] \label{def:Stieltjes}

\begin{itemize}
\item[(i)] Let  ${\mathfrak{f}}_z(x):=\frac{1}{x-z}$ ($x\in\R$, $z\in\C\setminus\R$).
\item[(ii)]
Let, for $z\in\C\setminus\R$, 
\BEQ U_t^N(z):=\langle X_t^N,{\mathfrak{f}}_z\rangle\equiv ({\cal S}X_t^N)(z)=\sum_{i=1}^N \frac{1}{\lambda^i_t-z}\EEQ
and 
\BEQ U_t(z):=\langle X_t,{\mathfrak{f}}_z\rangle\equiv ({\cal S}X_t)(z)=\int \frac{\rho_t(x)}{x-z}\, dx\EEQ
be the Stieltjes transform of $X_t^N$, resp. $X_t$.
\end{itemize}
\end{Definition}


\begin{Definition}[upper half-plane]
\begin{enumerate}
\item  Let $\Pi^+:=\{z\in\C\ |\ \Im(z)>0\}$.
\item For $b_{max}>0$, let $\Pi^+_{b_{max}}:=\{z\in \C\ |\ 0<\Im (z)<b_{max}\}$.
\item Let $\Pi^-:=-\Pi^+$, $\Pi^-_{b_{max}}:=-\Pi^+_{b_{max}}$ and $\Pi:=\Pi^+ \uplus \Pi^-,\ \Pi_{b_{max}}:=\Pi^+_{b_{max}}\uplus \Pi^-_{b_{max}}$.
\end{enumerate}
\end{Definition}

\begin{Definition} Let, for $p\in[1,+\infty]$ and $b_{max}>0$,
\BEQ L^p(\Pi_{b_{max}}):=\{h:\Pi_{b_{max}}\to\C\ |\  h(\bar{z})=\overline{h(z)}\ 
 (z\in \Pi^+_{b_{max}}) \, {\mathrm{\ and\ }} ||h||_{L^p(\Pi_{b_{max}})}<\infty\}, \EEQ
 where 
\BEQ ||h||_{L^p(\Pi_{b_{max}})}:= \left(\int_{-\infty}^{+\infty} da\,  \int_{-b_{max}}^{b_{max}} db\,  |h(a,b)|^p\right)^{1/p} \ (p<\infty), \qquad\!\!\!\!
||h||_{L^{\infty}(\Pi_{b_{max}})}:=\sup_{z\in\Pi_{b_{max}}} |h(z)|. \EEQ
\end{Definition}

\noindent The value of $b_{max}$ is unessential, so we fix some constant $b_{max}>0$ 
(e.g. $b_{max}=1$) and omit the $b_{max}$-dependence in the estimates.

\begin{Definition}[{\bf Stieltjes decomposition}] \label{def:Stieltjes-decomposition}
Let $\kappa=0,1,2,\ldots$ 
and $h\in L^1(\Pi_{b_{max}})$. Let $R=R(T)$ be as in Proposition \ref{prop:LD}. We say that $f:\R\to\R$ {\em has  Stieltjes decomposition $h$ of order $\kappa$ and cut-off $b_{max}$ on $[-R,R]$
} if, for all $|x|\le R$, 
\BEQ f(x)=({\cal C}^{\kappa}h)(x):=\int_{-\infty}^{+\infty} da \, \int_{-b_{max}}^{b_{max}} db\, (-\II b)\,  \frac{|b|^{\kappa}}{(1+\kappa)!}\,  \  
{\mathfrak{f}}_z(x) h(a,b). \label{eq:Stieltjes-decomposition}  \EEQ
\end{Definition}
Thanks to the symmetry condition, $h(\bar{z})=\overline{h(z)}$, 
(\ref{eq:Stieltjes-decomposition})  may be rewritten in the form
\BEQ ({\cal C}^{\kappa}h)(x)=2 \int_{-\infty}^{+\infty} da \, \int_{0}^{b_{max}}  db\, \frac{b^{1+\kappa}}{(1+\kappa)!}\,  \  
\Im\left[{\mathfrak{f}}_z(x) h(a,b)\right],   \EEQ
from which it is apparent that $f$ is indeed real-valued. 

\medskip\noindent Various Stieltjes decompositions, following Israelsson \cite{Isr},  have been constructed in \cite{Unt1}. The simplest one consists in defining 
 $h:(a,b)\mapsto  {\cal K}^{\kappa}_{b_{max}}(f)(a)$
where ${\cal K}^{\kappa}_{b_{max}} : f\mapsto {\cal F}^{-1}(K^{\kappa}_{b_{max}})\ast f$
is the Fourier multiplication operator by $K^{\kappa}_{b_{max}}(s):=\Big(
2\int_0^{b_{max}} db\, |b|^{1+\kappa} \ \cdot\ {\cal F}(\Im ({\mathfrak{f}}_{\II b}))(s)
\Big)^{-1}$. When $\kappa$ is {\em even}, it is proved (see \cite{Unt1}, (2.13)) that
\BEQ {\cal K}^{\kappa}_{b_{max}} = (1+\underline{\cal K}^{\kappa}_{b_{max}}) 
\Big( -(\kappa+1)! \, \partial_x^{2+\kappa} + (-1)^{\kappa/2} (2+\kappa) b_{max}^{-2-\kappa} \Big),  \label{eq:Kkappa}
\EEQ
where $|||\underline{\cal K}^{\kappa}_{b_{max}}|||_{L^1(\R)\to L^1(\R)},
|||\underline{\cal K}^{\kappa}_{b_{max}}|||_{L^{\infty}(\R)\to L^{\infty}(\R)}=O(1)$.
(We shall only need to consider
$\kappa=0$ in the present article). Since in the sequel we want to focus on narrow strips
around the real axis, one might think of taking the limit $b_{max}\to 0$. However, this
introduces awkward boundary terms. Instead we fix $b_{max}>0$ and define
$h:(a,b)\mapsto  e^{-b/\eps} {\cal K}^{\kappa}_{b_{max},\eps}$ ($\eps>0$), where ${\cal K}^{\kappa}_{b_{max},\eps}$ is the Fourier multiplication operator by 
$K^{\kappa}_{b_{max},\eps}(s):=\Big(
\int_0^{b_{max}} db\, b^{1+\kappa} e^{-b/\eps} \ \cdot\ {\cal F}(\Im ({\mathfrak{f}}_{\II b}))(s)
\Big)^{-1}$. Similarly to (\ref{eq:Kkappa}), we get (specifically for $\kappa=0$)
\BEQ {\cal K}^0_{b_{max},\eps} = (1+\underline{\cal K}^0_{b_{max},\eps}) 
\Big({\cal F}^{-1}((|s|+\eps^{-1})^2)\ast\Big)  \label{eq:Kkappaeps}
\EEQ
where $|||\underline{\cal K}^0_{b_{max},\eps}|||_{L^1(\R)\to L^1(\R)},
|||\underline{\cal K}^0_{b_{max},\eps}|||_{L^{\infty}(\R)\to L^{\infty}(\R)}=O(1)$.
The Fourier multiplication operator in the r.h.s. of (\ref{eq:Kkappaeps}) is not a differential operator any more: 
\BEQ \Big({\cal F}^{-1}((|s|+\eps^{-1})^2)\ast f\Big)(x)= \eps^{-2} \Big( \int_{-\infty}^{+\infty} f(y)\, dy \Big) + \Big(2\eps^{-1}
H\partial_x  -\partial_x^2 \Big) f(x)  \label{eq:Kkappaepsto0} \EEQ
where $H$ is the Hilbert transform (see Appendix). Note that the most singular term 
in $O(\eps^{-2})$ is simply a constant.

\bigskip\noindent In \cite{Unt1}, we wrote down an explicit 
time-dependent operator ${\cal H}(t)$   such that the right-hand side of
(\ref{eq:intro-PDE-f-asymptotic}) for $f_t$ decomposed as 
\BEQ  f_t(x)\equiv ({\cal C}^{\kappa}h_t)(x)=\int_{-\infty}^{+\infty} da  \, \int_{-b_{max}}^{b_{max}}  db\,  (-\II b)\,  \frac{|b|^{\kappa}}{(1+\kappa)!} \,  \  {\mathfrak{f}}_{z}(x) \, h_t(a,b) \label{eq:generator} \EEQ
(see Definition \ref{def:Stieltjes-decomposition}) is equal to  
\BEQ {\cal C}^{\kappa}({\cal H}(t)h_t)(x)=\int_{-\infty}^{+\infty} da\, \int_{-b_{max}}^{b_{max}}  db\, (-\II b)\,  \frac{|b|^{\kappa}}{(1+\kappa)!}\, 
 {\mathfrak{f}}_z(x)\,  {\cal H}(t)(h_t)(a,b) \EEQ

Note that, since
 Stieltjes decompositions are not unique, the operator ${\cal H}(t)$ is very
 under-determined. The essential features
 of the operator ${\cal H}(t)$ chosen in \cite{Unt1}  are recapitulated in Appendix, see section 4; in particular, for $\kappa\ge 0$, ${\cal H}(t)$ is the generator of a
 time-inhomogeneous semi-group of $L^p$, $p\ge 1$, which is a bounded perturbation of a
 transport operator.
Moving around the operator ${\cal H}(t)$ to the function ${\mathfrak{f}}_z(x)$, one
obtains an operator ${\cal L}(t)$ which is a "twisted adjoint" of ${\cal H}(t)$,

\BEQ {\cal L}(t):=w^{-1}(a,b) {\cal H}^{\dagger}(t) w(a,b) \label{eq:L}
\EEQ
 with $w(a,b):=(-\II b)
\frac{|b|^{\kappa}}{(1+\kappa)!}$ (see \cite{Unt1}, eq. (3.18) for details).   For
$\kappa=-1$ (at least formally), ${\cal L}(t)={\cal H}^{\dagger}(t)$, and $h_t$ may be
directly interpreted as a density in $\C\setminus\R$ exactly as in (\ref{eq:density}),
so that ${\cal L}(t)$ is the direct generalization of (\ref{intro:L}) to an arbitrary
potential.


\subsection{The stationary covariance kernel in the Hermite case}


We rewrite in appropriate coordinates the formulas for the covariance found by
Israelsson in the stationary regime, in the Hermite case, i.e.  when $V$ is harmonic $(V(x)=x^2/2)$ and $\beta=2$. 

The covariance kernel
$\Lambda$ of the Stieltjes transform
of the fluctuation field satisfies
the following obvious properties due to stationarity, where one assumes $t_1\ge t_2$,
\BEQ  \Lambda(t_1,z_1;t_2,z_2)=:
\Lambda(\Del t;z_1,z_2)=\lim_{t\to + \infty}
\Lambda(t+\Del t,z_1;t,z_2)\EEQ
with $\Del t:=t_1-t_2\ge 0$.

\medskip\noindent The covariance function $\Lambda$ can be found by taking the $t\to\infty$ limit in (Bender \cite{Ben}, Corollary 2.4 p. 7) with $\beta=2$, (take $\sigma=1/2$ and
$U(t,z)=\frac{1}{\sqrt{2}} U_{Bender}(t,\frac{z}{\sqrt{2}})$) 

\BEQ \Lambda(\Del t;z_1,z_2)=e^{-\Del t} \frac{\half f'_{\mu}(\frac{z_1}{\sqrt{2}})
\, \cdot\,  \half f'_{\mu}(\frac{z_2}{\sqrt{2}})}{2  \left(\half\,   (\frac{1}{\sqrt{2}}f_{\mu}(\frac{z_1}{\sqrt{2}})) \, \cdot\,  (\frac{1}{\sqrt{2}} f_{\mu}(\frac{z_2}{\sqrt{2}}))\,  e^{-\Del t}-1  \right)^2} \EEQ

where:

\begin{itemize}
\item[(i)] $\frac{1}{\sqrt{2}}f_{\mu}(\frac{z}{\sqrt{2}})=U_{eq}(z)=-z+\sqrt{z^2-2}$ is the Stieltjes transform of Wigner's
semi-circle law; the boundary values on the support $[-\sqrt{2},\sqrt{2}]$  are $U_{eq}(x\pm \II 0)=-x\pm\II\sqrt{2-x^2}$.
\item[(ii)] $\frac{1}{2}f'_{\mu}(\frac{z}{\sqrt{2}})=U'_{eq}(z)= \frac{z}{\sqrt{z^2-2}}-1
$ is its derivative, with boundary values on the support $U'_{eq}(x\pm\II 0)=\frac{x\mp\II \sqrt{2-x_1^2}}{\pm\II\sqrt{2-x_1^2}}$.
\end{itemize}

\medskip\noindent Letting $x_j\pm\II \sqrt{2-x_j^2}=: \sqrt{2}\,   e^{\pm i\theta_j}$  $(0\le \theta_j\le 
\pi)$ for convenience, i.e. $x_j=\sqrt{2} \cos(\theta_j)$, $\pi\rho_{eq}(x)=\sqrt{2-x_j^2}=\sqrt{2} \sin \theta_j$,
we get
\BEQ \Lambda(\Del t;x_1\pm\II 0,x_2\pm\II 0)=\frac{-1}{2\sin(\pm\theta_1)\, \sin(\pm\theta_2)} \frac{e^{-\II(\pm\theta_1\pm\theta_2)-\Del t}}{(e^{-\II(\pm\theta_1\pm\theta_2)-\Del t}-1)^2}, \label{eq:Lambda12} \EEQ
from which
\BEA &&  2\sin(\theta_1)\sin(\theta_2) g_{1,2}(t_1,x_1;t_2,x_2)=\frac{1}{4\pi^2}  \left[ \frac{e^{\II(\theta_1+\theta_2)-\Del t}}{(e^{\II(\theta_1+\theta_2)-\Del t}-1)^2}  +\frac{e^{\II(-\theta_1+\theta_2)-\Del t}}{(e^{\II(-\theta_1+\theta_2)-\Del t}-1)^2} +
\right. \nonumber\\
&& \qquad \left. 
+ \frac{e^{\II(\theta_1-\theta_2)-\Del t}}{(e^{\II(\theta_1-\theta_2)-\Del t}-1)^2} +
\frac{e^{-\II(\theta_1+\theta_2)-\Del t}}{(e^{-\II(\theta_1+\theta_2)-\Del t}-1)^2}
 \right] \nonumber\\
 && \qquad  =-\frac{1}{8\pi^2} \Re \left[ \frac{1}{\sin^2 \frac{\theta_1-\theta_2+\II\Del t}{2}} + \frac{1}{\sin^2 \frac{\theta_1+\theta_2+\II\Del t}{2}} \right]. \label{eq:g-Hermite} \EEA
  When $\Del t=0$, we find 
 \BEQ 2\sin(\theta_1) \sin(\theta_2) g_{1,2}(t,\theta_1;t,\theta_2)=-\frac{1}{8\pi^2} \left(\frac{1}{\sin^2(\frac{\theta_1-\theta_2}{2})} +  \frac{1}{\sin^2(\frac{\theta_1+\theta_2}{2})} \right) =  -\frac{1}{2\pi^2}
  \frac{1-\cos\theta_1\cos\theta_2}{(\cos\theta_1-\cos\theta_2)^2}  \label{eq:3.20}  \EEQ
Eq. (\ref{eq:3.20}) is in agreement with Johansson's formula for equilibrium fluctuations (compare with the kernel of the operator $h\mapsto \del^h$ of eq. (2.10) in \cite{Joh}, Theorem
2.4), in our case
\BEQ g_{1,2}(t,x_1;t,x_2)=\frac{1}{2\pi^2} \frac{1}{\sqrt{2-x_1^2}} \partial_{x_2}\left(
\frac{\sqrt{2-x_2^2}}{x_2-x_1}\right) =-\frac{1}{2\pi^2} \frac{1}{\sqrt{(2-x_1^2)(2-x_2^2)}} \frac{2-x_1 x_2}{(x_2-x_1)^2}. \EEQ

\vskip 2cm \noindent  We shall now be able to formulate a {\em hydrodynamic fluctuation equation}, compare with H. Spohn's formulas (3.30) or (4.14) in \cite{Spo}.

\medskip\noindent
From (\ref{eq:Lambda12}), one sees trivially  that  the boundary values $\Lambda(\Del t;x_1\pm\II 0, x_2+\pm\II 0)$ satisfy
the following evolution equation,
\BEQ (\partial_t\mp \II \partial_{\theta_1})\left(\sin(\theta_1)\, \Lambda(\Del t;x_1\pm\II 0,\cdot) \right)=0 \EEQ
or equivalently, considering the coordinate-transformed 
\BEQ \tilde{g}^{\eps_1,\eps_2}_{1,2}(\Del t;\theta_1, \theta_2):= 2\sin(\theta_1)\sin(\theta_2) g^{-\eps_1,-\eps_2}_{1,2}(\Del t;x_1,x_2)  \label{eq:gtilde}
\EEQ so that $\tilde{g}_{1,2}^{\eps_1,\eps_2}(\cdot;\theta_1,
\theta_2) d\theta_1\, d\theta_2= g^{-\eps_1,-\eps_2}_{1,2}(\cdot;x_1,x_2)
dx_1\, dx_2$,
\BEQ (\partial_{t_1}\pm\II \partial_{\theta_1}) \tilde{g}^{\pm,\cdot}_{1,2}(t_1,\theta_1
;\cdot) =0.   \label{eq:pre-hydro-Hermite} \EEQ
Since $\cos:(0,\pi)-\II\R_+^*\to \Pi_+$ and $\cos:(0,\pi)+\II \R_+^*\to \Pi_-$ are
biholomorphisms, the boundary value identity $\frac{1}{2\II\pi}(\tilde{g}^+ -\tilde{g}^-)
\equiv g$ on the real line may equally be seen as a boundary value identity 
 on the
unit circle in the variable $\tilde{z}_1:=e^{\II\theta_1}$, with $\tilde{g}^+$, resp.
$\tilde{g}^-$ extending on $\{|z_1|<1\}$, resp. $\{|z_1|>1\}$. 

\bigskip\noindent  Summing over $\eps_1,\eps_2$, we obtain the
{\em hydrodynamic fluctuation equation}
\BEQ \partial_{t_1} \tilde{g}_{1,2}(t_1,\theta_1;t_2,\theta_2)=
\partial_{\theta_1} (H_1\tilde{g}_{1,2})(t_1,\theta_1;t_2,\theta_2).
\label{eq:hydro-Hermite}
\EEQ
where $H_1$ is the {\em periodic Hilbert transform} acting on the first variable
(see Appendix, section 5). 
Since $\tilde{g}_{1,2}(t_1,\theta_1;t_2,\theta_2)$ is invariant under  the parity symmetries
$\theta_i\to -\theta_i$, $i=1,2$, we bother only about the first term
$\tilde{g}_{-}(t_1,\theta_1;t_2,\theta_2):=-\frac{1}{8\pi^2}\Re \frac{1}{\sin^2\frac{\theta_1-\theta_2+\II \Del t}{2}}$ in (\ref{eq:g-Hermite}),  and add the other term,
$\tilde{g}_{+}(t_1,\theta_1;t_2,\theta_2)=\tilde{g}_-(t_1,\theta_1;t_2,-\theta_2)$ 
by hand in the end. 
Since $\frac{\partial}{\partial\theta_1} \Big(\half\cot(\frac{\theta_1-
\theta_2}{2}) \Big)=-\frac{1}{4} \frac{1}{\sin^2(\frac{\theta_1-\theta_2}{2})}$,
compare with  (\ref{eq:Hper}),  we get in Fourier modes (see Appendix, section 6)
\BEQ \hat{K}_{\infty}(n)=\frac{1}{2\pi^2} |n|,  \qquad \frac{d}{d\Del t}\big|_{\Del t=0} \hat{K}_{\infty}(t+\Del t,t;n)=|n|
\hat{K}_{\infty}(n)
\qquad (n\in\Z) \EEQ
where $K_{\infty}(\cdot,\cdot)=\tilde{g}_{-}(t,\cdot;t,\cdot)$ is the stationary
covariance kernel. Thus {\em the asymptotic fluctuation process $Y$ is the stationary
solution of the  Ornstein-Uhlenbeck (linear Langevin) equation}
\BEQ \partial_t Y(t,\theta)=- \sqrt{-\partial^2/\partial\theta^2}\,   Y(t,\theta) + \frac{1}{\pi} \,  \partial_{\theta} \eta(t,\theta),  \label{eq:sqrt-D2} \EEQ
where $\sqrt{-\partial^2/\partial\theta^2}=\partial_{\theta}H$ is the convolution operator acting
by multiplication on 
Fourier modes, viz. $\widehat{\sqrt{-\partial^2/\partial\theta^2}} \hat{\phi}_n=|n|
\hat{\phi}_n$, and $\eta=\eta(t,\theta)$ is a $2\pi$-periodic space-time white noise 
admitting the parity symmetry, $\eta(t,-\theta)=\eta(t,\theta)$.

\vskip 1cm

\bigskip\noindent It is very instructive to compute the short-distance asymptotics in a scaled limit, $\Del t=\eps\del t_{12}\to 0$, $x_1-x_2=\eps\del x_{12}\to 0$. Formula (\ref{eq:g-Hermite}) implies in the angular coordinates
\BEQ \tilde{g}_{1,2}(t+\eps \del t_{12}, \theta+\eps\del \theta_{12};t,\theta)\sim_{\eps\to 0}
-\frac{1}{2\pi^2} \eps^{-2} \Re \Big[ \frac{1}{(\del\theta_{12}-\II \del t_{12})^2} \Big]  \EEQ
independently of $\theta$, from which
\BEA g_{1,2}(t+\eps \del t_{12},x+\eps \del x_{12};t,x)&\sim_{\eps\to 0}& -\frac{1}{2\pi^2}  \eps^{-2}  \Re\left[ \frac{1}{(\del x_{12}+\II
\sqrt{2-x^2}\, \del t_{12})^2}\right]  \nonumber\\
&=&-\frac{1}{2\pi^2} \eps^{-2}\Re \left[ \frac{1}{(\del x_{12}+\II\pi \rho_{eq}(x)\del t_{12})^2}\right]. \label{eq:short-distance} \EEA
Note that only the first term in the r.-h.s. of (\ref{eq:g-Hermite}) contributes to (\ref{eq:short-distance}).


\subsection{PDE for the covariance kernel: the general case}


We shall now derive a  PDE for $g_{1,2}^{+,\pm}$ in whole generality. (A PDE
for $g_{1,2}^{-,\pm}$ is then obtained by complex-conjugating the first
space coordinate.)

\begin{Theorem}[hydrodynamic fluctuation equation for general $V$] \label{th:g-PDE}
The kernel $g_{1,2}^{+,\pm}(t_1,x_1;t_2,x_2)$ satisfies the following PDE in a
weak sense,
\BEQ \partial_{t_1}g^{+,\pm}_{1,2}(t_1,x_1;t_2,x_2)=- \partial_{x_1} \Big( \Big(\frac{\beta}{2} U_{t_1}(x_1+\II 0)+V'(x_1)\Big) 
g^{+,\pm}_{1,2}(t_1,x_1;t_2,x_2) \Big), \label{eq:g-PDE} \EEQ
that is, for any smooth, compactly supported test function $\psi=\psi(x_1)$,
\BEQ \partial_{t_1} \int dx_1\, \psi(x_1)\, g_{1,2}^{+,\pm}(t_1,x_1;\cdot)=
\int dx_1\, \psi'(x_1)\, \Big(\frac{\beta}{2} U_{t_1}(x_1+\II 0)+V'(x_1)\Big) 
g^{+,\pm}_{1,2}(t_1,x_1;\cdot).  \label{eq:g-PDE-test}
\EEQ
\end{Theorem}

{\bf Remark.} The product $U_{t_1}(x_1+\II 0) g_{1,2}^{\pm,\pm}(t_1,x_1;\cdot)$ makes
sense as a distribution because both $x_1\mapsto U_{t_1}(x_1+\II 0)$ and 
$x_1\mapsto g_{1,2}^{\pm,\pm}(t_1,x_1;\cdot)$ are obtained by convolution with respect
to the function $x\mapsto \frac{1}{x+\II 0}$, hence have Fourier support
$\subset \supp \Big({\cal F}( x\mapsto \frac{1}{x+\II 0})\Big)=\R_+$.

{\bf Proof.} A short but non rigorous proof goes as follows.  Fix $\kappa=-1$. Since\\ 
$b_1({\cal H}_{nonlocal}^{0,-1}(t)h)(a_1,b_1)=O(b_1)\to_{b_1\to 0} 0$
(see (\ref{eq:H-dec})), we consider the limit when $N\to\infty$
and $b_1\to 0^+$  of the 
  characteristic equations associated to ${\cal L}_{transport}:=
  {\cal H}_{transport}^{\dagger}$, see (\ref{eq:L}) and below, thus obtaining directly
  the solution  of the evolution equation with terminal condition $f_{t_1}={\mathfrak{f}}_{z_1}$, where $z_1\equiv a_{t_1}+\II b_{t_1}$. One finds
 \BEA && {\cal L}_{transport}=\Big(v_{hor}(t,z)\partial_a+v_{vert}(t,z)\partial_b
 +\tau^{-1}(t,z)\Big)^{\dagger} \nonumber\\
 &&\qquad =-v_{hor}(t,z)\partial_a-v_{vert}(t,z)\partial_b+
 \Big(\tau^{-1}(t,z)-\frac{\partial v_{hor}}{\partial a}(t,z)-
 \frac{\partial v_{vert}}{\partial b}(t,z) \Big).  \nonumber\\ 
 \EEA   
 Explicit formulas
  (\ref{eq:da(tau)/dtau},\ref{eq:db(tau)/dtau},\ref{eq:dEsharp(tau)/dtau}) for
  $v_{hor},v_{vert},\tau^{-1}$ yield (as follows from easy explicit computation, of from
  \cite{Unt1}, eq. (3.15),(3.41), (3.45) and (3.48), where one has set $b\equiv 0^+$)

\BEQ \frac{da_t}{dt}=-\frac{\beta}{2}\Re U_t(a_t+\II 0)-V'(a_t) \EEQ
\BEQ \frac{db_t}{dt}=-\frac{\beta}{2}\Im U_t(a_t+\II 0) \EEQ
\BEQ \frac{dc_t}{dt}=\left[-\frac{\beta}{2} U'_t(a_t+\II 0)-V''(a_t) \right]c_t
\label{eq:g-char} \EEQ
 Consider these to be the characteristics of a generalized transport operator ${\cal L}_{hol}$ acting on a function ${\mathfrak{f}}_{z_1}$
{\em analytic} on $\Pi_+$, so that $\partial_{z_1}\equiv\partial_{z_1}+\partial_{\bar{z}_1}\equiv \partial_{x_1}$: then 
\BEQ {\cal L}_{hol}(t)=-\left(\frac{\beta}{2}U_t(x_1+\II 0)+V'(x_1)\right)\partial_{x_1} - \left(\frac{\beta}{2}U'_t(x_1+\II 0)+
V''(x_1)\right),
\label{eq:Lhol} \EEQ
which is exactly the operator featuring in the r.-h.s. of (\ref{eq:g-PDE}), acting on the
$x_1$-variable.
This makes it possible to keep $b_t\equiv 0^+$ during the time evolution.

\medskip\noindent Consider now the time-evolution of
\BEQ \phi(t_1,t_2;\lambda_1,\lambda_2;z_1,z_2)=:\esper\left[e^{\II \lambda_1\langle Y_{t_1},{\mathfrak{f}}_{z_1}\rangle}
\, e^{\II \lambda_2\langle Y_{t_2},{\mathfrak{f}}_{z_2}\rangle}\right] -
\esper\left[e^{\II \lambda_1\langle Y_{t_1},{\mathfrak{f}}_{z_1}\rangle} \right]
\, \esper\left[ e^{\II \lambda_2\langle Y_{t_2},{\mathfrak{f}}_{z_2}\rangle}\right], \EEQ
with  $z_1=x_1+\II 0$ and $t_1\ge t_2$. 
Taylor expanding around $(\lambda_1,\lambda_2)=(0,0)$ yields
\BEQ \phi(t_1,t_2;\lambda_1,\lambda_2;z_1,z_2)=1-\lambda_1 \lambda_2 \Lambda(t_1,z_1;t_2,z_2)+O((|\lambda_1|+|\lambda_2|)^3). \EEQ
On the other hand, by (\ref{eq:intro-Gaussian}),
\BEA && \phi(t_1,t_2;\lambda_1,\lambda_2;z_1,z_2)= \exp\Big( \half \int_{t_2}^{t_1}  ds\, \Big( \II
(1-\frac{\beta}{2}) \langle X_s,\lambda_1 (f_s)''\rangle -\langle X_s,
(\lambda_1 (f_s)')^2 \rangle  \Big) \Big)\ \cdot \nonumber\\
&&\qquad \cdot\ \Big(  \esper\left[e^{\II \lambda_1\langle Y_{t_2},f_{t_2}\rangle}
\, e^{\II \lambda_2\langle Y_{t_2},{\mathfrak{f}}_{z_2}\rangle}\right] -
\esper\left[e^{\II \lambda_1\langle Y_{t_2},f_{t_2}\rangle} \right]  \, \esper\left[ e^{\II \lambda_2\langle Y_{t_2},{\mathfrak{f}}_{z_2}\rangle}\right]
\Big)  \label{eq:lambda1lambda2}
\EEA
where $f_s$ is the solution at time $s\le t_1$ of (\ref{eq:intro-PDE-f-asymptotic}) with terminal
condition $f_{t_1}\equiv {\mathfrak{f}}_{z_1}$.  The second line of 
(\ref{eq:lambda1lambda2}) is of the form  $-C(t_1,z_1;t_2,z_2)\lambda_1\lambda_2 +
O((|\lambda_1|+|\lambda_2|)^3)$, hence (by identification)  $C(t_1,z_1;t_2,z_2)\equiv \Lambda(t_1,z_1;t_2,z_2)$.
 Thus the time-evolution
of $\Lambda$ may be computed by considering the sole contribution of (\ref{eq:Lhol})
acting on the $x_1$ variable.

\bigskip\noindent For a rigorous proof we proceed instead as follows. Consider an
arbitrary terminal condition $f=f_{t_1}\in L^1(\R)\cap C^{\infty}(\R)$ such that
$\int_{-\infty}^{+\infty} f(y)\, dy=0$ (this can always be achieved by modifying $f$
outside of the support of the density), and rewrite
it as $f\equiv f_+ + f_-\equiv {\cal C}^0 (h)+{\cal C}^0 (\bar{h})$, where $h(z)\equiv h_{t_1}(z)={\bf 1}_{b>0} e^{-b/\eps} ({\cal K}_{1,\eps}^0 f)(z)$ (see (\ref{eq:Kkappaeps}),
(\ref{eq:Kkappaepsto0})) and
$\bar{h}(\bar{z}):=\overline{h(z)}$,  for $\eps>0$. Later on we let $\eps\to 0$ to obtain the
time-evolution  of $\Lambda$ near the real axis. Define
\BEQ  \Phi_f(t_1,t_2;\lambda_1,\lambda_2;z_2)=:\esper\left[e^{\II \lambda_1\langle Y_{t_1},f_+\rangle}
\, e^{\II \lambda_2\langle Y_{t_2},{\mathfrak{f}}_{z_2}\rangle}\right] -
\esper\left[e^{\II \lambda_1\langle Y_{t_1},f_+\rangle} \right]
\, \esper\left[ e^{\II \lambda_2\langle Y_{t_2},{\mathfrak{f}}_{z_2}\rangle}\right]. \EEQ
Taylor expanding around $(\lambda_1,\lambda_2)=(0,0)$ yields
\BEQ \Phi_f(t_1,t_2;\lambda_1,\lambda_2;z_2)=-\lambda_1 \lambda_2 \Lambda_{f}(t_1;t_2,z_2)+O((|\lambda_1|+|\lambda_2|)^3) 
\EEQ
where 
\BEQ  \Lambda_{f}(t_1;t_2,z_2):=\int da_1 \int_0^{b_{max}} (-\II b_1) e^{-|b_1|/\eps}
({\cal K}^0_{1,\eps}f)(a_1) \Lambda(t_1,z_1;t_2,z_2).
\EEQ
Eq.   (\ref{eq:lambda1lambda2}) for $\phi(t_1,t_2;\lambda_1,\lambda_2;z_1,z_2)$ holds
equally well for $\Phi_f(t_1,t_2;\lambda_1,\lambda_2;z_2)$, with only the terminal condition $f_{t_1}$ changing. Thus we get two expressions for the time-derivative
of $\Lambda_f$,
\BEA  \frac{\partial}{\partial t_1} \Lambda_f(t_1;t_2,z_2) &=& \int da_1 \int_0^{b_{max}} db_1 \, (-\II b_1)\,  e^{-b_1/\eps} ({\cal K}^0_{1,\eps}f)(a_1) \partial_{t_1}\Lambda(t_1,z_1;t_2,z_2)
\nonumber\\
&=&\int da_1 \int_0^{b_{max}} db_1 \, (-\II b_1)\,   {\cal H}^0(t_1) \Big(e^{-b_1/\eps} ({\cal K}^0_{1,\eps}f)(a_1) \Big) \Lambda(t_1,z_1;t_2,z_2). 
\label{eq:HKLambda} \nonumber\\
\EEA
The condition $\int_{-\infty}^{+\infty} f=0$ ensures that $({\cal K}^0_{1,\eps}f)(a_1)
\sim_{\eps\to 0} 2\eps^{-1}(Hf)'(a_1)$. 

\medskip\noindent Main terms in (\ref{eq:HKLambda}) are due to ${\cal H}^0_{transport}(t_1)\equiv v_{hor}(t,z_1)\partial_{a_1}+v_{vert}(t,z_1)\partial_{b_1}+\tau^0(t,z_1)$.
Consider first the  velocity terms; since $U_t\Big|_{\Pi_+}$ extends continuously
to the real axis by hypothesis, $e^{-|b_1|/\eps} v_{hor}(t,z_1)=
e^{-|b_1|/\eps} v^{asympt}_{hor}(t,a_1)+ o(1)$  when 
$\eps\to 0$, where $v^{asympt}_{hor}(t,a_1):=\frac{\beta}{2} \Re U_t(a_1+\II 0)+V'(a_1)$.
Then (integrating by parts) the adjoint operator $(v_{vert}(t,z_1)\partial_{b_1})^{\dagger}= -v_{vert}(t,z_1)\partial_{b_1}-\frac{\partial v_{vert}(t,z_1)}{\partial b_1}$, $(\cdots)$, 
acts on the product $b_1 \Lambda(t_1,z_1;t_2,z_2)$. Since $z_1\mapsto \Lambda(t_1,z_1;\cdot)$ is holomorphic, $\partial_{b_1}\Lambda(t_1,z_1;\cdot)\equiv \II  \partial_{a_1}
\Lambda(t_1,z_1;\cdot)$. Thus the action of $v_{vert}(t,z_1)\partial_{b_1}=\Big((v_{vert}(t,z_1)\partial_{b_1})^{\dagger}
\Big)^{\dagger}$ is equivalent to that of  the
transport operator $\II v_{vert}(t,z_1)\partial_{a_1}$. 
Now $e^{-|b_1|/\eps} (v_{hor}(t,z_1)+\II v_{vert}(t,z_1))=e^{-|b_1|/\eps} v^{asympt}(t,a_1)+o(1)$, with $v^{asympt}(t,a_1)=\frac{\beta}{2}U_t(a_1+\II 0)+V'(a_1)$, a function
whose product with $\Lambda(t_1,z_1;t_2,z_2)$ is well-defined (see Remark before the
proof). 
Then:
\BEA &&\int_0^{b_{max}} db_1\, (-\II b_1)\,   e^{-b_1/\eps} \int da_1 \, v^{asympt}(t,a_1)
\Lambda(t_1,z_1;t_2,z_2)  \partial_{a_1}  ({\cal K}^0_{1,\eps}f)(a_1) \nonumber\\
&&=\int_0^{b_{max}} db_1\, (-\II b_1)\, e^{-b_1/\eps} \esper\Big[ ({\cal S}Y)(t_2,z_2)\,   \langle Y_{t_1},x\mapsto \int da_1\, 
v^{asympt}(t,a_1) \frac{ \partial_{a_1}  ({\cal K}^0_{1,\eps}f)(a_1)}{x-a_1-\II b_1} \rangle  \Big] \nonumber\\
&&\sim_{\eps\to 0} -\II\eps^{-1} I(\eps) 
 \int da_1\,  \psi'(a_1)\Big(\frac{\beta}{2} U_{t_1}(x_1+\II 0)+V'(x_1)\Big) 
\Lambda(t_1,x_1+\II 0;t_2,z_2)  \nonumber\\  \label{eq:cov-velocity}
\EEA
where $\psi(a_1):= 2(Hf)'(a_1)$ and $I(\eps):= \int_0^{b_{max}} db_1 \, b_1 e^{-b_1/\eps}\sim_{\eps\to 0} \eps^2$ is simply a coefficient. Missing multiplicative terms are easily checked (as we
did in the above "short" proof) to
compensate all terms in  (\ref{eq:dEsharp(tau)/dtau}) with $\kappa=0$, except the last
one, $\frac{dc_t^0}{dt}=\cdots - \II V'''(a_t)b_t c_t^0$; however, the corresponding multiplication
operator $-\II V'''(a_1)b_1$, as well as $b_1 {\cal H}^{1,0}_{nonloca}(t_1)$,   vanish 
in the limit $\eps\to 0$ when multiplied by $e^{-|b_1|/\eps}$.

\medskip\noindent On the other hand, working directly on the first line of eq. (\ref{eq:HKLambda}) yields
 $\frac{\partial}{\partial t_1}\Lambda_f(t_1;t_2,z_2)\sim_{\eps\to 0} -\II\eps^{-1}  I(\eps) \int da_1\,  \psi(a_1) \partial_{t_1}\Lambda(t_1,a_1+\II 0;t_2,z_2)$. Comparing the latter expression with (\ref{eq:cov-velocity}) yields (\ref{eq:g-PDE-test}).
\hfill \eop

\vskip 1cm
 
\bigskip\noindent {\em Assume from now on that $V$ extends analytically to
an entire function $V:\C\to\C$} (e.g. $V$ is a polynomial). Then  the above equation (\ref{eq:g-PDE}) may be solved formally in terms
of its initial condition, namely,  the one-time covariance kernel $g_{1,2}^{+,\pm}(t_2,\cdot;t_2,\cdot)$, as
follows. First one uses the adjoint equation in order to obtain a transport operator
without multiplicative term. Namely, letting $\psi_t$ solve the PDE
$\partial_t\psi_t(x)=-(\frac{\beta}{2}U_t(x+\II 0)+V'(x))\psi_t(x)$, we have
$\partial_{t_1} \langle \psi_{t_1},g^{+,\pm}(t_1,\cdot;t_2,x_2) \rangle=0$.  {\em 
Assume} $\psi_{t_1}(x)\equiv \psi_{t_1}(x+\II 0)$ is the boundary value of a 
function $\psi_{t_1}=\psi_{t_1}(z)$ holomorphic on $\Pi_+$; then, for $t\le t_1$,
$\psi_t(x)\equiv \psi_t(x+\II 0)$ is the boundary value of $\psi_t=\psi_t(z)$ solution of 
\BEQ \partial_t \psi_t(z)=-(\frac{\beta}{2} U_t(z)+V'(z))\psi_t(z) \label{eq:psi}.
\EEQ
 Denote
by $\Phi_t^{t_1}(z_1)$ $(0<t<t_1)$ the solution of the characteristic equation
on $\Pi_+$, 
$\dot{z}=-(\frac{\beta}{2} U_t(z)+V'(z))$ with terminal condition
 $z(t_1)\equiv z_1\in\Pi_+$.
Note that we considered exactly the same characteristics in \S 2.2 {\bf B}, where
 the notation $Z_{t_1}$ for $\Phi_t^{t_1}$ was used.
Then $\psi_{t_2}(z_1)=\psi_{t_1}(\Phi_{0}^{t_1}(z_1))$. Hence (considering the
scalar product $\langle \, \cdot \, , \, \cdot\,  \rangle$ in $L^2(\R)$)
\BEQ \langle \psi_{t_1}, g_{1,2}^{+,\pm}(t_1,\cdot;t_2,x_2) \rangle
 =\langle
\psi_{t_2}, 
g_{1,2}^{+,\pm}(t_2,\cdot;t_2,x_2)\rangle  \label{eq:g-PDE-solution}
\EEQ

\medskip
\noindent Eq. (\ref{eq:g-PDE-solution}) does not suffice to characterize the law of the
fluctuation process in general. However, Theorem \ref{th:g-PDE} can easily be
extended to provide a full answer. Namely, letting $\eps_i=\pm$, $i=1,2$
\begin{itemize}
\item[(i)] $\partial_t g_{1,2}^{\eps_1,\eps_2}(t,x_1;t,x_2)= 
\Big({\cal L}_{\eps_1}^{(1)}(t)+
{\cal L}_{\eps_2}^{(2)}(t) \Big) g_{1,2}^{\eps_1,\eps_2}(t,x_1;t,x_2)$, where
\BEQ {\cal L}_{+}^{(i)}(t)=-\left(\frac{\beta}{2}U_t(x_i+\II 0)+V'(x_i)\right)\partial_{x_i} - \left(\frac{\beta}{2}U'_t(x_i+\II 0)+
V''(x_i)\right).
\label{eq:Lplus} \EEQ
is ${\cal L}_{hol}$, see (\ref{eq:Lhol}), acting on the $x_i$-variable, while
${\cal L}_-:=\overline{{\cal L}_+}$. Solving as in (\ref{eq:g-PDE-solution}) by
following
characteristics on both space variables simultaneously, one obtains \\ 
$g_{1,2}^{\eps_1,\eps_2}(t,x_1;t,x_2)$ in terms of the initial covariance $g_{1,2}^{\eps_1,\eps_2}(0,\cdot;
0,\cdot)$. 

\item[(ii)] The one-point function $\esper[Y(x)]$ does not follow from the above
computations; actually,
Johansson  gave a general but not very explicit formula for $\esper[Y(x)]$
at equilibrium, showing that it  vanishes for $\beta=2$ but not for $\beta\not=2$ in
general (see \cite{Joh}, Theorem 2.4, Remark 2.5 and eq. (3.54)). Following the method used in the proof of Theorem \ref{th:g-PDE}, one can consider
the time-evolution of the generating function
\BEQ \phi(t;\lambda;z):= \esper\left[ e^{\II\lambda \langle Y_t,{\mathfrak{f}}_z\rangle}
\right]
\EEQ
with $z=x+\II 0$. Taylor expanding around $\lambda=0$ yields $\phi(t;\lambda;z)=1+\II
\lambda \esper[({\cal S}Y_t)(z)] + O(\lambda^2)$. On the other hand,
\BEQ \phi(t;\lambda;z)=\exp\Big( \half \int_0^{t}  ds\, \Big( \II
(1-\frac{\beta}{2}) \langle X_s,\lambda (f_s)''\rangle -\langle X_s,
(\lambda (f_s)')^2 \rangle  \Big) \Big) \cdot\   \esper\left[e^{\II \lambda\langle Y_0,f_0\rangle}
\right]  \label{eq:lambda}
\EEQ
where $f_s$ is the solution at time $s\le t$ of (\ref{eq:intro-PDE-f-asymptotic}) with terminal
condition $f_{t}\equiv {\mathfrak{f}}_{z}$. Differentiating w.r. to $t$ and Taylor
expanding to order 1 in $\lambda$ yields
$\frac{d}{dt} \esper[({\cal S}Y_t)(z)]= \Big( \frac{\II}{2} (1-\frac{\beta}{2})
 U''_t(z) + {\cal L}_{\pm}(t,z)  \Big) \esper[({\cal S}Y_t)(z)]$,
where ${\cal L}_{\pm}={\cal L}_+$ if $z\in\Pi_+$, resp. ${\cal L}_-$ if $z\in\Pi_-$,
a generalized transport equation on $\C\setminus\R$ with the same characteristics as
above.
\end{itemize}

\bigskip
\noindent How explicit can these formulas be made ? One  may of course try to answer
this question through
case-by-case inspection. Let us point out at two specific but sufficiently general
cases. The first one is the {\em harmonic case}, i.e. $V(x)=\frac{x^2}{2}$,
treated in an exhaustive way by Bender \cite{Ben} (see in particular Theorem 2.3) for an arbitrary parameter
$\beta>1$ and an arbitrary initial condition. Though the mapping $\Phi_t^{t_1}$
is explicit (see (\ref{eq:Z})), the {\em inverse mapping}, $(\Phi_t^{t_1})^{-1}$, of
course, is not in general. It requires some skill to provide explicit formulas
not relying on the use of $(\Phi_t^{t_1})^{-1}$, see e.g. the beautiful
result using Schwartzian derivatives (\cite{Ben}, Theorem 2.7) for $\Cov( \langle Y_{t_1},F_1\rangle, \langle Y_{t_2},
F_2\rangle)$ when $F_1,F_2$ are bounded analytic functions on a neighbourhood of the
real axis.  The second one is the {\em stationary case}, where $\beta$ and $V$ are
general but $\rho_t=\rho_{eq}$ is assumed to be the equilibrium measure. This is the
subject of the next subsection.


\subsection{Solution of the PDE in the stationary case}


We restrict to the {\em stationary case} in this subsection, and assume as stated
before that $V$ {\em extends analytically to an entire function} $V:\C\to\C$. Let us first
state two essential facts.  First,
the {\em universality} (up to simple scaling and translation) {\em of Johansson's formula for  equilibrium fluctuations}
implies, assuming that $\supp(\rho_{eq})=[-A,A]$ $(A>0)$: 

\BEA  g_{1,2}^{+,\pm}(t_1,x_1;t_1,x_2) &=& (\frac{\sqrt{2}}{A})^2 \Lambda_{(\ref{eq:Lambda12})}(0;\frac{\sqrt{2}}{A}(x_1+\II 0),\frac{\sqrt{2}}{A}(x_2\pm \II 0)) \nonumber\\
&=& (\frac{\sqrt{2}}{ A})^2 \frac{1}{8\sin(\theta_1)\sin(\pm \theta_2) \sin^2 \frac{\theta_1\pm\theta_2}{2}}
\EEA
where $\Lambda_{(\ref{eq:Lambda12})}$ is as in (\ref{eq:Lambda12}), and   $A\cos(\theta_j)=x_j, A\sin(\theta_j)=\sqrt{A^2-x_j^2}$, $j=1,2$ is up to 
scaling the change of variables used
in the Hermite case.   Second,  using (\ref{eq:cut-eq1}) and (\ref{eq:H}), (\ref{eq:pm}),
\BEQ U_{eq}(x\pm \II 0)=-\frac{2}{\beta}V'(x)\pm\II \pi \rho_{eq}(x),\EEQ
which may also be interpreted as saying that $\rho_{eq}$ extends analytically 
to $\Pi_{\pm}$ as $\frac{1}{\pm\II\pi}(U_{eq}(z)+\frac{2}{\beta}V'(z))$.
Therefore, {\em Theorem \ref{th:g-PDE} may be restated in this simple form,}
\BEQ \partial_{t_1}  g_{1,2}^{+,\pm}(t_1,x_1;t_2,x_2)=-\II \pi \frac{\beta}{2} \partial_{x_1} \Big(\rho_{eq}(x_1) g_{1,2}^{+,\pm}(t_1,x_1;t_2,x_2) \Big).
\label{eq:g-PDE-stat}
\EEQ
which  generalizes (\ref{eq:hydro-Hermite}).
Solving (\ref{eq:g-PDE-stat}) for short time and $x_1\to x_2$, with initial
condition $(t_1=t_2)$ given  by
  Johansson's equilibrium formula, one finds the same short-distance
asymptotics as in (\ref{eq:short-distance}), namely,
\BEQ g_{1,2}^{+,+}(t+\eps \del t_{12},x+\eps \del x_{12};t,x)\sim_{\eps\to 0} -\frac{1}{4\pi^2}  \eps^{-2}  \left[ \frac{1}{(\del x_{12}+\II\pi \rho_{eq}(x)\del t_{12})^2}\right]. \label{eq:short-distance-rho} \EEQ
See discussion in the Introduction.

\bigskip\noindent The  hydrodynamic fluctuation
equation  (\ref{eq:g-PDE-stat}) may be solved as follows:

\begin{Theorem}[equilibrium fluctuation covariance kernel] \label{th:g-stat}
Let $A\cos(\theta_1)=x_1$, $A\sin(\theta_1)=\sqrt{2-x_1^2}$ as above, and:
\begin{itemize}
\item[(i)] $G=G(z)$ be the
analytic continuation to $\Pi_+$ of the function $G$ defined on the support of $\rho_{eq}$ by $G(x):=\frac{2}{\beta} \int_0^x\frac{dy}{\rho_{eq}(y)};$
\item[(ii)] $F_0^{\pm}(\cdot,x_2):=F^{\pm}_0(z,x_2)$ be the analytic continuation to $\Pi_+$ of the function 
$F^{\pm}_0(x_1,x_2):=\rho_{eq}(x_1) g_{1,2}^{+,\pm}(0,x_1;0,x_2)= (\frac{\sqrt{2}}{A})^2  \rho_{eq}(A\cos\theta_1) \frac{1}{8\sin(\theta_1)
\sin(\pm\theta_2)\sin^2 \frac{\theta_1\pm \theta_2}{2}}.$
\end{itemize}

Then
\BEQ g_{1,2}^{+,\pm}(t,x_1;0,x_2)=\frac{1}{\rho_{eq}(x_1)} F^{\pm}_0(G^{-1}(G(x_1)+\II \pi t),
x_2).  \label{eq:g-stat} \EEQ

\end{Theorem}

\noindent{\bf Remark.} Extend continuously $G$, resp. $F^{\pm}_0$ to $[-A,A]\uplus \Pi^+\uplus
\Pi^-$ by letting $G(\bar{z}):=\overline{G(z)}$, resp. $F^{\pm}_0(\bar{z}):=
\overline{F^{\pm}_0(z)}$. Then
$$ g_{1,2}^{-,\pm}(t,x_1;0,x_2)=\overline{g_{1,2}^{+,\pm}(t,x_1;0,x_2)}= \frac{1}{\rho_{eq}(x_1)} F^{\pm}_0(G^{-1}(G(x_1)-\II \pi t),
x_2)$$
satisfies the time-reversed evolution equation $\overline{(\ref{eq:g-PDE-stat})}$,
$$ \partial_{t_1}  g_{1,2}^{-,\pm}(t_1,x_1;t_2,x_2)=+\II \pi \frac{\beta}{2} \partial_{x_1} \Big(\rho_{eq}(x_1) g_{1,2}^{-,\pm}(t_1,x_1;t_2,x_2) \Big).$$

\medskip\noindent
{\bf Proof.} The product $\psi(t,x_1;0,x_2):=\rho_{eq}(x_1)g_{1,2}^{+,\pm}(t,x_1;0,x_2)$ solves the transport equation
$\partial_{t}\psi=-\II\pi  \frac{\beta}{2} \rho_{eq}(x_1) \partial_{x_1}\psi$ with initial condition
$\psi(0,x_1;0,x_2)=\rho_{eq}(x_1)g_{1,2}^{+,\pm}(0,x_1;0,x_2)=F^{\pm}_0(x_1,x_2)$. Since
$-\II\pi\frac{\beta}{2}\rho(x_1)=\lim_{z\to x_1,z\in\Pi_+} -(\frac{\beta}{2}U_{eq}(z)+V'(z))$, $\psi(t,x_1;0,x_2)$ 
may be obtained exactly as in (\ref{eq:psi}) by considering  the 
time-homogeneous characteristic equation in $\Pi_+$, $\dot{z}_1=-(\frac{\beta}{2}U_{eq}(z)+V'(z))$ with terminal condition $z_1(t)\equiv x_1$.  Solving by quadrature, one gets
\BEQ t=-\int_{z_1(0)}^{x_1} \frac{dw}{\frac{\beta}{2}U_{eq}(w)+V'(w)} \equiv 
\frac{1}{\II\pi}(G(z_1(0))-G(x_1)).\EEQ  
Whence $\psi(t,x_1;0,x_2)=F^{\pm}_0(z_1(0),x_2)$, equivalent to (\ref{eq:g-PDE-stat}).
 \hfill \eop

\bigskip\noindent Using (\ref{eq:gggg}), one obtains:

\BEA && g_{1,2}(t_1,x_1;t_2,x_2)=-\frac{1}{2\pi^2} \frac{1}{\rho_{eq}(x_1)} \ \Re \Big[
F_0^+(G^{-1}(G(x_1)+\II\pi (t_1-t_2)),x_2) \nonumber\\
&&\qquad\qquad\qquad\qquad\qquad  - F_0^-(G^{-1}(G(x_1)+\II\pi (t_1-t_2)),x_2)  \Big]
\EEA

\bigskip\noindent Using time-translation invariance for $g,g^{\pm,\pm}$ and
reversibility for $g$, a more symmetric-looking formula can be obtained. Namely,
assume $t:=t_1-t_2>0$ and let $t'\in(\frac{t}{2},t)$. Then
\BEA g_{1,2}(t_2+t',x_1;t_2,x_2) &=& g_{1,2}(t'-\frac{t}{2},x_1;-\frac{t}{2},x_2)=
g_{1,2}(-\frac{t}{2},x_1;t'-\frac{t}{2},x_2) \nonumber\\
&=& -\frac{1}{2\pi^2} \Re [g^{+,+}_{1,2}-g^{+,-}_{1,2}]( -\frac{t}{2},x_1;t'-\frac{t}{2},x_2).
\EEA
Then $\psi'(-\frac{t}{2},x_1;t'-\frac{t}{2},x_2):=\rho_{eq}(x_2)
g^{+,\pm}(-\frac{t}{2},x_1;t'-\frac{t}{2},x_2)$ solves  the transport equation  
 $\partial_{t'}\psi'=-\II \pi \frac{\beta}{2} \rho_{eq}(x_2) \partial_{x_2} \psi'$
 with initial condition $\psi'(-\frac{t}{2},x_1;0,x_2)=\rho_{eq}(x_2) 
 \overline{g^{+,\pm}(\frac{t}{2},x_1;0,x_2)}$ by the previous Remark, solved as in
 the  proof of Theorem \ref{th:g-stat} by solving a characteristic equation
 in the {\em second} space variable $x_2$. The conjugation plays no r\^ole in the end upon taking the real part. The sign in the denominator $\sin(\pm\theta_2)=\pm 
 \sin\theta_2$ of $F_0^{\pm}$ can also be removed since (\ref{eq:gggg}) is an
 alternate sum. Hence the result of Theorem \ref{th:g-stat} can be reformulated
 as follows:
 
\begin{Corollary} \label{cor:g-stat}
Let 
\BEQ \tilde{g}^{\pm}(x_1,x_2):=(\frac{\sqrt{2}}{A})^2  \rho_{eq}(A\cos\theta_1)\rho_{eq}(A\cos\theta_2) \frac{1}{8\sin(\theta_1)
\sin(\theta_2)\sin^2 \frac{\theta_1\pm \theta_2}{2}}.\EEQ
Then
\BEA && g_{1,2}(t_1,x_1;t_2,x_2)=-\frac{1}{2\pi^2} \frac{1}{\rho_{eq}(x_1)\rho_{eq}(x_2)} \ \Re \Big[ \tilde{g}^+\Big(G^{-1}(G(x_1)+\II \frac{\pi}{2} (t_1-t_2)), \nonumber\\
&& \ \  G^{-1}(G(x_2)+\II \frac{\pi}{2} (t_1-t_2))\Big) \  +\  \tilde{g}^-\Big(G^{-1}(G(x_1)+\II \frac{\pi}{2} (t_1-t_2)), G^{-1}(G(x_2)-\II \frac{\pi}{2} (t_1-t_2))\Big)  \Big].
\nonumber\\
\EEA
\end{Corollary}

The notation $\tilde{g}^{\pm}_{1,2}$ for the covariance multiplied by the
density is reminiscent of (\ref{eq:gtilde}). In the Hermite case (see \S 3.2), 
the $\rho_{eq}$-factors in the denominator of $\tilde{g}^{\pm}$ cancel the
sines in the denominator, $G(x)=\pi\theta$ and $G^{-1}(G(x)+\II\pi\frac{t}{2})=
G^{-1}(\pi(\theta+\II \frac{t}{2}))=\sqrt{2} \cos(\theta+\II\frac{t}{2})$, thus
time-evolution amounts to an imaginary translation of the angle coordinates
$\theta_1,\theta_2$, and one easily retrieves (\ref{eq:g-Hermite}).



\vskip 1cm

\bigskip\noindent Let us illustrate this explicit formula in the case when
$V(x)=\frac{1}{4} t^4+ \frac{c}{2} t^2+d$ is quartic. The additive constant $d$ does not play any r\^ole, so we may forget it. The equilibrium density $\rho_{eq}$ and its Stieltjes
transform $U_{eq}$, when $\supp(\rho_{eq})$ is {\em connected}, are given by a general explicit integral formula, see e.g. \cite{Joh}, eq. (3.9), which can be solved when
$V$ is polynomial. The result in the particular case when  $V$ is quartic can be
found e.g. in Johansson \cite{Joh}, Example 3.2; note the slight discrepancy of notations
in that article with respect to ours, namely, $U_{Johansson}=-U$, $V_{Johansson}=\frac{4}{\beta} V$, $T_{Johansson}=\frac{4}{\beta} T$.  Thus, letting $\supp(\rho_{eq})=:
[-A,A]$, $A=A(c)$, we find:
\BEQ U_{eq}(z)=\frac{2}{\beta} \Big(-V'(z)+(z^2+\half A^2+c)\sqrt{z^2-A^2}\Big),\  
\rho_{eq}(x)=\frac{2}{\beta\pi}(x^2+\half A^2+c)\sqrt{A^2-x^2}\, {\bf 1}_{|x|<A}.\EEQ 
This is the same density as Johansson's up to a rescaling,\\ $\rho_{eq}(x)=(\frac{\beta}{2})^{-1/4} \rho_{eq,Johansson}((\frac{\beta}{2})^{-1/4}x; (\frac{\beta}{2})^{-1/4}A,
(\frac{\beta}{2})^{-1/2}c)$, which means that\\
$A(c)=(\frac{\beta}{2})^{1/4} A_{Johansson}((\frac{\beta}{2})^{1/2}c)=
(\frac{\beta}{2})^{1/2} \sqrt{-\frac{2}{3}c+2\sqrt{\frac{1}{9}c^2+\frac{8}{3\beta}}}.$
Then (reversing the arrow of time for commodity of notation) we  solve  forward characteristics
with fixed initial condition,
\BEQ \dot{z}=+(\frac{\beta}{2}U(z)+V'(z))=(z^2+\half A^2+c)\sqrt{z^2-A^2}, \qquad 
z(0)\equiv x_1\in\Pi_+.\EEQ
Changing variables, $\theta\equiv \arcsin (x/A)$, then $h:=\cotan \theta=\sqrt{\frac{A^2}{z^2}-1}$, one finds by quadrature
\BEQ t=\II \tau \Big\{ \arctan \Big(C\sqrt{(A^2/z_t^2)-1}\Big)  - \arctan\Big(C\sqrt{(A^2/x_1^2)-1}\Big)\Big\} \EEQ
where $\tau:=\frac{1}{\sqrt{(\frac{3}{2}A^2+c)(\half A^2+c)}},\ 
 C:=\sqrt{\frac{\half A^2+c}{\frac{3}{2}A^2+c}}$. In the notations of
 Theorem \ref{th:g-stat} or Corollary \ref{cor:g-stat}, one has found:
\BEQ G(z)=\pi\tau \arctan \Big(C\sqrt{(A^2/z^2)-1}\Big). \EEQ
 Inverting this formula, we get
\BEA  z_t&\equiv & G^{-1}(G(x_1+\II\pi t)) \nonumber\\
& =&  A \Big( 1+ C^{-2} \tan^2(-\II \frac{t}{\tau}+\arctan(C\sqrt{(A^2/x_1^2)-1})
\Big)^{-1/2} \nonumber\\
 & =&  A \Big( 1+ C^{-2} \Big[ \frac{-\II \tanh(\frac{t}{\tau})+C\sqrt{(A^2/x_1^2)-1}}{1+
\II C \tanh(\frac{t}{\tau})\sqrt{(A^2/x_1^2)-1}} \Big]^2
\Big)^{-1/2} \nonumber\\
&=& A\Big(1+
\II C \tanh(\frac{t}{\tau})\sqrt{(A^2/x_1^2)-1}\Big) \ \cdot  \ \Big[ \ (1-C^2 \tanh^2(t/\tau))
\frac{A^2}{x_1^2} + \tanh^2(t/\tau) (C^2-C^{-2})  \nonumber\\
&& \qquad\qquad + 2\II (C-C^{-1})\tanh(t/\tau)
\sqrt{(A^2/x_1^2)-1} \ \   \Big]^{-1/2}.    \label{eq:g-PDE-stat-quartic}
\EEA


\section{Appendix. Generator and semi-group estimates}


A large part of the work in our previous article \cite{Unt1} has been to write down
explicitly a time-dependent  operator ${\cal H}^{\kappa}(t)$ (called: {\em 
generator}) such that, assuming 
$f_T={\cal C}^{\kappa}(h_T)$, the function $f_t={\cal C}^{\kappa}(h_t)$ with $h_t$
solution for $t\le T$ of the backwards evolution equation
\BEQ \frac{dh}{dt}(t;a,b)={\cal H}^{\kappa}(t)(h(t))(a,b), \qquad h(T)\equiv h_T \EEQ
is solution of (\ref{eq:intro-PDE-f-asymptotic}).

\medskip\noindent The operator ${\cal H}^{\kappa}(t)$ exhibited in \cite{Unt1} is
of the following form:
\BEQ {\cal H}^{\kappa}(t)(h(t))(a,b)\equiv {\cal H}^{\kappa}_{transport}(t)(h(t))(a,b)+
b {\cal H}^{\kappa+1,\kappa}_{nonlocal}(h(t))(a,b),   \label{eq:H-dec} 
\EEQ
where (assuming $\supp (\rho_t)\subset [-R,R]$, $0\le t\le T$): 

\begin{enumerate}
\item
${\cal H}^{\kappa}_{transport}$
 is a (generalized)  transport
operator (see \cite{Unt1}, section 6 for a brief exposition of the characteristic method,
and \S 3.6 for the formulas below, where we have taken the limit $N\to\infty$),

\BEQ {\cal H}^{\kappa}_{transport}(t)=v_{hor}(t,z)\frac{\partial}{\partial a} +
v_{vert}(t,z) \frac{\partial}{\partial b}+\tau^{\kappa}(t,z) \EEQ
with  associated
characteristics on $[-3R,3R]\times[-b_{max},b_{max}]$,
\BEQ \frac{da_t}{dt}=v_{hor}(t,z_t)\equiv \frac{\beta}{2} \Re(U_{t}(a_t+\II b_t)) +V'(a_t)-\half V'''(a_t) b_t^2   \label{eq:da(tau)/dtau} \EEQ
\BEQ \frac{db_t}{dt}=v_{vert}(t,z_t)\equiv \frac{\beta}{2} \Im(U_{t}(a_t+\II b_t))+ V''(a_t)b_t \label{eq:db(tau)/dtau} \EEQ 
\BEA &&  \frac{dc^{\kappa}_{t}}{dt}=\tau^{\kappa}(t,z_t)\equiv \Big[\frac{\beta}{2} \Big(
\frac{1+\kappa}{b_t}\Im(U_{t}(a_t+\II b_t)) +  (\bar{U}_{t})'(a_t+\II b_t)  \Big)  \nonumber\\
&& \qquad\qquad\qquad  +(2+\kappa) V''(a_t) 
-\II V'''(a_t) b_t\Big] c^{\kappa}_{t}. \label{eq:dEsharp(tau)/dtau} \EEA

Let $(a_t,b_t)$ be the solution of (\ref{eq:da(tau)/dtau},\ref{eq:db(tau)/dtau})
with terminal condition $(a_T,b_T)\equiv (a,b)$, and $c_t^{\kappa}\equiv \exp\Big( -\int_t^T
\tau^{\kappa}(a_s,b_s)\, ds\Big)$ the solution of (\ref{eq:dEsharp(tau)/dtau}) with terminal
condition $c_T^{\kappa}\equiv 1$. Then the solution of
the evolution equation $\frac{\partial f_t}{\partial t}(x)={\cal H}^{\kappa}_{transport}(t) f_t(x)$ with terminal condition $f_T\equiv f$ is: $f_t(a,b)=c_t^{\kappa} f(a_t,b_t)$.

\item (non-local term)
 \BEQ   |||{\cal H}_{nonlocal}^{\kappa+1;\kappa}|||_{
(L^1\cap L^{\infty})(\Pi_{b_{max}}) \to (L^1\cap L^{\infty})(\Pi_{b_{max}})} = O( ||V'||_{8+\kappa,[-3R,3R]} ),\EEQ
 (see \cite{Unt1}, eq. (4.26),(4.27)).

\end{enumerate}

\bigskip
\noindent The transport operator ${\cal H}_{transport}^{\kappa}(t)$ can be exponentiated backward
in time for $\kappa\ge 0$ as results from the sign of $\Re \tau^{\kappa}$: namely,

\begin{Proposition}(see \cite{Unt1}, Lemma 3.5) \label{lem:U-transport}
Let $u_T\in (L^1\cap L^{\infty})([-3R,3R]\times(0,b_{max}])$ and $\kappa=0,1,2\ldots$ 
Then the backward evolution equation $\frac{du}{dt}={\cal H}^{\kappa}_{transport}(t)u(t), 
u\big|_{t=T}=u_T$ $(0\le t\le T)$ has a unique solution $u(t):=U^{\kappa}_{transport}(t,T)u_T$, such that 
\BEQ ||u_t||_{L^1([-3R,3R]\times(0,b_{max}])}\le  ||u_T||_{L^1([-3R,3R]\times(0,b_{max}])} \EEQ
\BEQ ||u_t||_{L^{\infty}([-3R,3R]\times(0,b_{max}])}\le  ||u_T||_{L^{\infty}([-3R,3R]\times(0,b_{max}])} \EEQ
\end{Proposition}


\section{Appendix. Stieltjes and Hilbert transforms} \label{sec:Stieltjes-Hilbert}


We collect in this section some definitions and elementary properties concerning
Stieltjes  and Hilbert transforms, in the periodic and in the non-periodic cases.


\subsection{Non-periodic case}


 We make use of the Fourier transform normalized as follows, ${\cal F}(f)(s)=\int_{-\infty}^{+\infty} f(x) e^{-\II xs}\, dx$. 

\medskip

\noindent Let, for $z=a+\II b\in\C\setminus\R$,  ${\mathfrak{f}}_z(x)=\frac{1}{x-z}
 (x\in\R)$, so that the {\em Stieltjes transform} of $\phi:\R\to\R$ is $({\cal S}\phi)(z):=\langle \phi,{\mathfrak{f}}_z\rangle\equiv \int dx\, \phi(x){\mathfrak{f}}_z(x)$. 
Note that $\Im({\mathfrak{f}}_z)(x)=\frac{b}{(x-a)^2+b^2}$ is of the same sign as $b$.
The Plemelj
formula, $\frac{1}{x-\II 0}=p.v.\left(\frac{1}{x}\right)+\II \pi \del_0$
implies the  following boundary values,
\BEQ \frac{1}{2\II\pi} (\phi_+(x)-\phi_-(x))= \frac{1}{\pi}\Im \phi_+(x)=
 \phi(x)  \label{eq:bv-} \EEQ
with $\phi_{\pm}(x):=\lim_{b\to 0^+} ({\cal S}\phi)(x\pm\II b)$, and
\BEQ \frac{1}{2\pi} (\phi_+(x)+\phi_-(x))=-(H\phi)(x), \label{eq:bv+}\EEQ
where $H\phi$ is the {\em Hilbert transform} of $\phi$,
\BEQ (H\phi)(x):= \frac{1}{\pi} p.v.\ \int \frac{dy}{x-y} \phi(y). \label{eq:H} 
\EEQ
 Conversely,
\BEQ \frac{1}{\pi} \phi_{\pm}(x)=-(H\phi)(x)\pm\II \phi(x). \label{eq:pm} \EEQ
Since
$${\cal F}{\mathfrak{f}}_{\II b}(s)=2\II\pi \sgn(b) e^{-|b|\, |s|} {\bf 1}_{\sgn(s)=-\sgn(b)},$$ 
the Fourier kernel of the Hilbert transform is ${\cal F}H(s)=-\II \sgn(s)$,
implying in particular $H^2=-I$. Applying this identity  to a function $\phi$ yields
\BEQ ({\cal F}\phi^+)(s)=2\II\pi {\bf 1}_{s<0} ({\cal F}\phi)(s), \qquad 
 ({\cal F}\phi^-)(s)=-2\II\pi {\bf 1}_{s>0} ({\cal F}\phi)(s).
\EEQ

\medskip\noindent {\em An essential property of the Stieltjes transform
$({\cal S}\rho)(z):=\langle\rho,{\mathfrak{f}}_z\rangle$ of a density is the following},
\BEQ  \Im ({\cal S}\rho)(z)=\langle \rho,x\mapsto \frac{b}{(x-a)^2+b^2}\rangle >0,
\qquad z\in\Pi_+   \label{eq:positivity} \EEQ
from which it follows that the functions $U_t={\cal S}\rho_t$ map $\Pi_+$ into
$\Pi_+$ and $\Pi_-$ into $\Pi_-$.


\subsection{Periodic case}


The Fourier integral is replaced by Fourier series,  $\phi(\theta)=\sum_{n\in\Z} c_n(\phi)
e^{\II n\theta}$ with 
\BEQ c_n(\phi)\equiv \hat{\phi}(n):=\frac{1}{2\pi} \int_0^{2\pi} d\theta\, \phi(\theta) e^{-\II n\theta}.
\EEQ 
{\em We consider only functions $\phi$ with vanishing means, i.e. such that
$c_0(\phi)=\frac{1}{2\pi}\int_0^{2\pi} \phi\equiv 0.$} The Stieltjes transform (still
denoted ${\cal S}$) is now
a Cauchy integral on the unit circle, $({\cal S}\phi)(z):=\oint_{|\zeta|=1} d\zeta\, \frac{\phi(\zeta)}{\zeta-z}$ ($|z|\not=1$).  In terms of the Fourier coefficients,
\BEQ ({\cal S}\phi)(z)=\sum_{n\ge 1} c_n(\phi) z^n \qquad  (|z|<1),  \qquad -\sum_{n\le -1} c_n(\phi) z^n \qquad  (|z|>1).
\EEQ 
Letting $\phi_{\pm}(\theta):=\lim_{r\to 1^{\mp}} ({\cal S}\phi)(re^{\II\theta})$, one
has: 
$$\phi_+(\theta)=2\II\pi\sum_{n\ge 1} c_n(\phi) e^{\II n\theta}, \qquad  \phi_-(\theta)=-2\II\pi\sum_{n\le  -1} c_n(\phi) e^{\II n\theta}.$$
 By analogy with the real line case,
we let 
\BEQ (H\phi)(\theta):=-\frac{1}{2\pi}(\phi_+(\theta)+\phi_-(\theta))=\sum_{n\ge 1}
(-\II) c_n(\phi)e^{\II n\theta} + \sum_{n\le -1} (+\II) c_n(\phi)e^{\II n\theta}, \EEQ
also given, using $\cot \frac{\theta}{2}=2 \Re \left\{ (-\II)\sum_{n\ge 1} e^{\II n\theta}\right\}$, by the following convolution kernel, 
\BEQ (H\phi)(\theta)=\frac{1}{2\pi} p.v. \int_0^{2\pi} dt\, f(t) \cot(\frac{\theta-t}{2}),\EEQ
and get $\phi=\frac{1}{2\II\pi} (\phi_+-\phi_-), H^2=-I$.
Letting  $\ccotan(\theta):=\half\cot \frac{\theta}{2},$
we obtain an equivalent formula,
\BEQ (H\phi)(\theta)=\frac{1}{\pi} p.v. \int_0^{2\pi} dt\, f(t) \ccotan(\theta-t),
\label{eq:Hper}
\EEQ
making it plain that the periodic Hilbert transformation is a rational generalization of
the Hilbert transform (\ref{eq:H}) on the real line.



\section{On Ornstein-Uhlenbeck processes}


An Ornstein-Uhlenbeck process is a (Hilbert-space valued) stochastic process
$Y(t)$ satisfying a linear stochastic differential equation of the form
\BEQ \dot{Y}(t)=-AY(t)+\Sigma \, \eta(t)  \label{eq:OU-SDE} \EEQ
where $\eta$ is delta-correlated white noise, time-derivative of a Wiener process, 
and $A,\Sigma$ are some operators; see \cite{Hai}, \S 5 for details. If $Y:\R_+\to\R$
is one-dimensional, and $\Sigma=\sqrt{T}>0$,  $Y(t)$ modelizes either  the velocity of a massive Brownian particle under
the influence of friction, or the position of an infinitely massive Brownian particle
submitted to friction and to a harmonic potential $V(Y)=\half AY^2$; in the first
interpretation, $A$ is the friction coefficient. In both cases $T$ plays the r\^ole
of  a temperature,
as appears in the Maxwell-like equilibrium distribution $e^{-AY^2/2T}=e^{-V(Y)/T}$.
In our context $Y=Y(t,x)$ is the random fluctuation process, $\eta=\eta(t,y)$ is
space-time white noise, and (\ref{eq:OU-SDE}) is a Langevin equation for $Y$.
Under adequate assumptions, notably on the analytic properties and long-time
behavior of the  semi-group $e^{-tA}$, $t\ge 0$ generated
by $A$, this equation has a unique stationary measure $\mu_{\infty}$, and the law $\mu_t$ of
$Y_t$ converges to $\mu_{\infty}$ for any reasonable initial measure  $\mu_0$. Furthermore, $\mu_{\infty}$ is Gaussian, with covariance kernel $K_{\infty}=K_{\infty}(x,y)$ defined
uniquely by 
\BEQ \mbox{Sym} (K_{\infty} A^{\dagger})=\half \Sigma\Sigma^{\dagger} \label{eq:OU-K} \EEQ
with $\mbox{Sym}(B):=\half(B+B^{\dagger})$ (see \cite{Hai}, Theorem 5.22). If $\Sigma,A$ are self-adjoint and commute, and $A\ge 0$, then (starting
 from any initial measure) $Y(t)=e^{-tA}Y(0)+e^{-tA} \int_0^t e^{sA} \Sigma\eta(s)$, so
 $K_{\infty}=\lim_{t\to \infty} \int_0^t ds\,  e^{-(t-s)A} \Sigma \Sigma^{\dagger} e^{-(t-s)A}=
 \half \Sigma^2/A$, confirming (\ref{eq:OU-K}).
 
  {\em Assume conversely} that some stationary Gaussian process 
 $Y(t)$ is given, with known two-time covariance kernel $K_{\infty}(t_1,x_1;t_2,x_2)=
 K_{\infty}(t_1-t_2;x_1,x_2)$. Then $Y$ is the solution of (\ref{eq:OU-SDE}) with
\BEQ AK_{\infty}=\frac{d}{dt}\big|_{t=0} K_{\infty}(t,0), \qquad \half \Sigma \Sigma^{\dagger}=\mbox{Sym}(K_{\infty} A^{\dagger}).\EEQ



\end{document}